\documentclass{amsart}
\usepackage{amssymb, amsmath, amscd, amsthm}
\usepackage[all]{xy}

       {\arraycolsep 0.14em\begin{eqnarray}}{\end{eqnarray}}
\newenvironment{Eqnarray*}%
       {\arraycolsep 0.14em\begin{eqnarray*}}{\end{eqnarray*}}
       {\arraycolsep 0.14em\begin{array}}{\end{array}}

\theoremstyle{plain}
\newtheorem{Thm}{Theorem}[section]
\newtheorem{Lem}[Thm]{Lemma}
\newtheorem{Prop}[Thm]{Proposition}

\newtheorem{Conj}[Thm]{Conjecture}

\theoremstyle{definition}
\newtheorem{Def}[Thm]{Definition}

\theoremstyle{remark}

\numberwithin{equation}{section}

\def\action{\, {\scriptscriptstyle \stackrel{\circ}{{}}} \, }

  \def\C{\mathbb C}  \def\R{\mathbb R}
\def\Z{\mathbb Z}  \def\H{\mathbb H}  \def\bP{\mathbb P}
 \def\bD{\mathbb D}

\def\a{\alpha}  \def\b{\beta}        
\def\c{\theta}  \def\l{\lambda}        \def\s{\sigma}
\def\t{\tau}    \def\vp{\varphi}       \def\y{\eta}
\def\z{\zeta}   \def\w{\omega}    \def\e{\varepsilon}
\def\x{\xi}       \def\k{\kappa}

        \def\D{\Delta}

   \def\cD{\mathcal D}   \def\cF{\mathcal F}
      
\def\cL{\mathcal L}       
\def\cO{\mathcal O}   \def\cZ{\mathcal Z}    \def\cH{\mathcal H}
\def\cP{\mathcal P}       \def\cS{\mathcal S}

      \def\Gm{\mathfrak m}
\def\Gq{\mathfrak q}      

\def\Gn{\mathfrak n}   \def\Gz{\mathfrak z}  \def\Gp{\mathfrak p}

\def\Imag{\operatorname{Im}}    \def\Real{\operatorname{Re}}

\def\RP{\R\bP}
\def\CP{\C\bP} 

\def\del{\partial}

\def\({ \left( }     \def\){ \right) }
\def\[{ \left[ }     \def\]{ \right] }
\def\^{ \wedge }
\def\<{ \left\langle } \def\>{ \right\rangle }
\def\hasira{\rule{0mm}{2eX}}

\def\pv{{\bf pv.}\!\!}
\def\chalf{\frac{\c}{2}}

\newcommand{\Frac}[2]{\frac{\textstyle #1}{\textstyle #2}}
\newcommand{\pd}[1]{\Frac{\del}{\del #1}}
\newcommand{\PD}[2]{\Frac{\del #1}{\del #2}}

\renewcommand{\bar}{\overline}

\newcounter{mynum}

\title{Singular self-dual Zollfrei Metrics and Twistor Correspondence}
\author{Fuminori Nakata}
\address{Graduate School of Mathematical Sciences, 
University of Tokyo, 3-8-1 Komaba Meguro, 
Tokyo 153-8914, Japan}
\email{nakata@ms.u-tokyo.ac.jp}
\thanks{The author's research is supported by Research Fellowship of 
 the Japan Society for the Promotion of Science for Young Scientists}

\begin{document}
\maketitle
\begin{abstract}
  We construct examples of singular self-dual Zollfrei metrics explicitly, 
  by patching a pair of Petean's self-dual split-signature metrics. 
  We prove that there is a natural one-to-one correspondence 
  between these singular metrics and a certain set of embeddings of $\RP^3$ to 
  $\CP^3$ which has one singular point. 
  This embedding corresponds to an odd function on $\R$
  that is rapidly decreasing and pure imaginary valued. 
  The one-to-one correspondence is explicitly given by using the Radon transform.
\end{abstract}

\tableofcontents

\section{Introduction}
A {\it Zollfrei metric}, which was introduced by V.~Guillemin (\cite{bib:Guill89}), 
is an indefinite metric of a manifold whose maximal null geodesics are all closed. 
C.~LeBrun and L.~J.~Mason investigated the self-dual Zollfrei metric of 
signature (2,2),  
and constructed its twistor correspondence (\cite{bib:LM05}). 
They proved that only $S^2\times S^2$ and $(S^2\times S^2)/\Z_2$ admit a 
self-dual Zollfrei conformal structure. 
Using the twistor correspondence, they also proved that
such structure on $(S^2\times S^2)/\Z_2$ is rigid, 
and, in contrast, $S^2\times S^2$ admits many of such structures. 
In the case of $S^2\times S^2$, the corresponding twistor space is given by a 
pair $(\CP^3,P)$, where $P$ is the image of a totally real embedding 
$\RP^3\rightarrow\CP^3$. 
Their theorem says that there is a one-to-one correspondence 
between self-dual Zollfrei conformal structures on $S^2\times S^2$ 
and the pairs $(\CP^3,P)$ at least in the neighborhoods of the standard structures. 

On the other hand, K.~P.~Tod (\cite{bib:Tod}) and H.~Kamada 
(\cite{bib:Kamada})
independently constructed infinitely many examples of 
$S^1$-invariant scalar-flat indefinite  
K\"{a}hler metrics on $S^2\times S^2$, which are automatically self-dual.  
Because Zollfrei condition is an open condition in the space of self-dual metrics 
(\cite{bib:LM05}), these examples should contain many self-dual Zollfrei metrics, 
and a natural problem here is to decide whether all of them are Zollfrei or not. 
These examples are explicitly written in closed form, so it might be possible 
to write down explicitly the twistor correspondence for such metrics. 
We are not going to pursue these questions in this article, 
instead of that, 
we generalize the formulation to admit certain singularity, 
and we construct explicit examples of singular self-dual Zollfrei metric, 
whose twistor correspondence is explicitly written down. 

While LeBrun and Mason's theorem stands only in the neighborhood of 
the standard metric 
because they are using inverse function theorem for Banach space, 
our examples contain many metrics which are far from the standard one. 
We use the Radon transform to write down the 
twistor correspondence for our singular metric. 

In \cite{bib:Petean}, J.~Petean classified the compact complex surfaces which 
admit indefinite K\"{a}hler-Einstein metrics. 
Petean constructed many self-dual metrics on $\R^4$ to show that the complex 
tori or the primary Kodaira surfaces admit many of such metrics. 
Our examples of the singular self-dual Zollfrei metric are constructed 
by patching a pair of J.~Petean's metrics on $\R^4$. 

Our main theorem is to establish the twistor correspondence 
for some class of singular metrics. 
We construct the singular metrics and the singular twistor spaces 
respecting certain fiber bundle structure over $S^2$ and $\CP^2$ respectively. 
The main theorem is proved respecting 
a fiber bundle structure over a natural double fibration 
$S^2\leftarrow Z \rightarrow \CP^2$ for some $Z$, 
which is the twistor correspondence for the standard Zoll metric on $S^2$.  
Remark that a Zoll metric on a smooth manifold is a Riemannian metric 
whose geodesics are all closed, 
and the simplest one is the standard metric on $S^2$. 
The general case of twistor correspondence for Zoll structure is established by 
LeBrun and Mason (\cite{bib:LM02}). 

The organization of the paper is the following. 
In Section \ref{Section:Standard}, 
we recall the statement of the twistor correspondence 
for Zoll or Zollfrei structure respectively following LeBrun and Mason 
(\cite{bib:LM02},\cite{bib:LM05}). 
We describe the twistor correspondence 
for the standard structure on $S^2$ or $S^2\times S^2$ 
by introducing local coordinates which we use later. 
In Section \ref{Section:Def_sing}, we introduce the definition of 
the singular self-dual Zollfrei metrics and the singular twistor spaces, and
we formulate a conjecture of the twistor correspondence between them 
(Conjecture \ref{Conj:twistor_corr}). 

In Section \ref{Section:sing_met}, we construct explicit examples of 
the singular self-dual Zollfrei metric by patching two Petean's metrics. 
Each example corresponds to an element of the set 
$\cS(\R^2)^{\text{sym}}$ defined below. 
Let $\cS(\R^2)$ the set of rapidly decreasing real functions on $\R^2$, 
and $\cS(\R^2)^{\text{sym}}$ be the subset of $\cS(\R^2)$ 
consisting of $SO(2)$-invariant elements, which we call axisymmetric functions.

In Sections \ref{Section:twistor} and \ref{Section:holo_disks}, 
we construct explicit examples of singular twistor space. 
Each twistor space corresponds to an element of $i\cS(\R)^{\text{odd}}$, 
where $i\cS(\R)^{\text{odd}}$ is the set of 
odd functions that are rapidly decreasing and pure imaginary valued. 
Our main theorem (Theorem \ref{Thm:twistor_correspondence}) says 
that our conjecture holds when we restrict to the above mentioned examples. 
This correspondence is explicitly given as a transform between 
$f(x)\in\cS(\R^2)^{\text{sym}}$ and $h(t)\in\cS(\R)^{\text{odd}}$, 
by using the Radon transform. 
In Section \ref{Section:corresp} we give the proof of 
Theorem \ref{Thm:twistor_correspondence}, and the appendix 
(Section \ref{Appendix}) is the review of the Radon transform.

\section{Standard model} \label{Section:Standard}

{\it Zoll projective structures :}
A Zoll metric on a smooth manifold $M$ is a Riemannian metric 
whose geodesics are all closed (cf.\cite{bib:Guill76}). 
An example of such a metric is the standard metric on $S^2$. 
A Zoll projective structure on $M$ is a projectively equivalence class 
of torsion-free connections on the tangent bundle $TM$ 
whose geodesics are all closed, where two torsion-free connections 
are said to be projectively equivalent if and only if 
they have exactly the same unparameterized geodesics (\cite{bib:LM02}). 
Each Zoll metric defines a Zoll projective structure. 

C.~LeBrun and L.~J.~Mason proved the following property in \cite{bib:LM02}; 
there is a natural one-to-one correspondence between 
\begin{itemize}
 \item equivalence classes of Zoll projective structures on $S^2$, 
 \item equivalence classes of totally real embeddings 
    $\iota: \RP^2 \rightarrow\CP^2 $, 
\end{itemize}
in neighborhoods of the standard projective structure on $S^2$ 
and the standard embedding of $\RP^2$. 
We call this correspondence the twistor correspondence for Zoll projective 
structures. 
This correspondence is characterized by the following condition; 
there is a double fibration $S^2 \overset{\Gp} \longleftarrow \cD 
\overset{\Gq}\longrightarrow \CP^2$ such that 
\begin{enumerate}
 \item  $\Gq$ is a continuous surjection 
             and $\Gp$ is a complex disk bundle, i.e. for each $x\in S^2$, 
			 $\cD_x=\Gp^{-1}(x)$ is biholomorphic to the complex unit disk, 
 \item  $\Gq_x : \cD_x \rightarrow\CP^2$ is holomorphic on the interior of 
   $\cD_x$, and $\Gq_x(\del \cD_x)\subset N$, 
   where $\Gq_x$ is the restriction of $\Gq$ on $\cD_x$ and $N=\iota(\RP^2)$, 
 \item the restriction of $\Gq$ on $\cD-\del\cD$ is bijective onto $\CP^2-N$, 
 \item $\{\Gp(\Gq^{-1}(y))\}_{y\in N}$ is equal to the set of geodesics 
   on $S^2$. 
\end{enumerate}
The conditions (2) and (3) say that $\{\Gq(\cD_x)\}_{x\in S^2}$ 
is a family of holomorphic disks on $\CP^2$ with boundaries on $N$ which foliate 
$\CP^2-N$. 

For the standard projective structure on $S^2$, 
the twistor correspondence is described by the following diagrams 
which are explained below: 
\begin{equation} \label{diag:double_fibration}
  \xymatrix{
    & \bP(TS^2) \ar[dl]_p \ar[dr]^q & & 
	& \bD(TS^2) \ar[dl]_\Gp \ar[dr]^\Gq & \\
    S^2 & & \RP^2 &	S^2 & & \CP^2 }
\end{equation}
We call the left diagram the real twistor correspondence 
and the right one the complex twistor correspondence. 

We denote $S^2=\{t\in\R^3 : \| t\| ^2=1\}$, 
$TS^2=\{(t,v)\in S^2\times\R^3 : \< t,v\> =0 \}$ 
and $\bP(TS^2)=\{(t,[v])\in S^2\times\RP^2 : \< t,v\> =0 \}$. 
Let $p:\bP(TS^2)\rightarrow S^2$ be the projection and 
$q(t,[v])=[t\times v]$, then we have the left diagram in 
(\ref{diag:double_fibration}). 

Here we introduce a local coordinate system. 
We take an open covering $\{D_+,D_-,W\}$ of $S^2$ where 
$D_\pm=\{t\in S^2 : \pm t_3>0 \}$ and $W=S^2-\{(0,0,\pm 1)\}$. 
We define local coordinates 
$(x^\pm_1,x^\pm_2)\in\R^2\simeq D_\pm$ 
and $(\a,\b)\in\R/2\pi\Z\times[-\frac{\pi}{2},\frac{\pi}{2}]\simeq W$ by 
$$ \begin{pmatrix} t_1 \\ t_2 \\ t_3 \end{pmatrix}  
  = \pm \( 1+(x^\pm_1)^2+(x^\pm_2)^2\)
     \begin{pmatrix} x^\pm_1 \\ x^\pm_2 \\ 1 \end{pmatrix} 
  = \begin{pmatrix} \cos\a\cos\b \\ \sin\a\cos\b \\ \sin\b \end{pmatrix}. $$
The circle bundle $\bP(TS^2)$ is trivial over $D_\pm$ and $W$, and 
we define the fiber coordinate $\z^\pm$ and $\x$ by 
\begin{equation} 
 \begin{array}{c}
 D_\pm\times(\R\cup\{\infty\}) \ni \(x_1^\pm,x_2^\pm, \z^\pm \) 
    \longleftrightarrow 
    \[ -\z^\pm\pd{x_1^\pm}+\pd{x_2^\pm} \] 
    \in \bP(TS^2)|_{D_\pm}, \\[2ex]
 W \times (\R\cup\{\infty\}) \ \ni \(\a,\b,\x\) 
    \longleftrightarrow
    \[ -\x\pd{\a}+\cos^2\b\, \pd{\b} \]
    \in \bP(TS^2)|_W. 
 \end{array}
\end{equation}
The coordinate change is given by 
\begin{equation} \label{coordinate_change}
 \begin{cases} 
   x_1^\pm = \cos\a\cot\b, \\
   x_2^\pm = \sin\a\cot\b,
 \end{cases} \qquad
 \z^\pm =\frac{-\x\sin\a\tan\b +\cos\a}{\x\cos\a\tan\b-\sin\a}. 
\end{equation}
In terms of these coordinates, the map $q:\bP(TS^2)\rightarrow\RP^2$ 
is described by 
\begin{equation}\label{map_q}
 \begin{aligned}
  \(x_1^\pm,x_2^\pm,\z^\pm\) & \longmapsto 
   \[ 1 : \z^\pm : -x_1^\pm -x_2^\pm \z^\pm \], \\
  \(\a,\b,\x\) \ & \longmapsto 
  \[\x\cos\a\tan\b-\sin\a : -\x\sin\a\tan\b +\cos\a : \x \]. 
 \end{aligned}
\end{equation} 

We define a complex disk bundle $\bD(TS^2)$ over $S^2$ as a closure of one 
of the two connected components of $\bP(T_\C S^2)-\bP(TS^2)$, 
where $\bP(T_\C S^2)$ is the complex projectivization of 
$T_\C S^2=TS^2\otimes\C$. 
The choice of the component is not essential because 
these components are canonically 
isomorphic by the complex conjugation. 

Extending real parameters $\z^\pm$ and $\x$ to the complex parameter 
inte upper or lower half plane $\H_\pm=\{z\in\C : \pm\Imag z \geq 0\}$, 
we can introduce the trivialization of $\bD(TS^2)$ by 
$$ \begin{aligned} 
 \bD(TS^2)|_{D_\pm} &\simeq D_\pm\times (\H_\pm \cup\{\infty\})
    \ni (x_1^\pm, x_2^\pm, \z^\pm), \\
 \bD(TS^2)|_W &\simeq W\times (\H_+ \cup\{\infty\})
    \ni (\a, \b, \x). \end{aligned} $$
The coordinate change is given by the same formulas as (\ref{coordinate_change}) 
with complex coordinates. 
Then the map $\Gq : \bD(TS^2) \rightarrow \CP^2$ is obtained by the 
analytic continuation of $q$, i.e. $\Gq$ is given by the same formula as 
(\ref{map_q}). 
It is easy to check that the above conditions (1), (2), (3) and (4) hold 
if we put $\cD=\bD(TS^2)$.  \\

{\it Self-dual Zollfrei metrics :}
A Zollfrei metric on a smooth manifold $M$ is an indefinite metric 
whose null geodesics are all closed (cf.\cite{bib:Guill89}). 
In \cite{bib:LM05}, LeBrun and Mason investigated self-dual Zollfrei neutral metrics
on four dimensional manifolds, 
where a neutral metric is an indefinite metric with signature $(n,n)$ which is 
also called an indefinite metric with split signature. 
An example of such a metric is the standard metric $g_0$ on 
$S^2\times S^2$ given by $g_0=\pi_1^*h_{S^2}-\pi_2^*h_{S^2}$, 
where $\pi_i$ is the projection to the $i$-th $S^2$ and $h_{S^2}$ is the 
standard Riemannian metric on $S^2$. 

LeBrun and Mason proved the following property in \cite{bib:LM05}; 
let $M=S^2\times S^2$ and $g$ be a self-dual Zollfrei neutral metric on $M$, 
then every {\it $\b$-surface} on $M$ is homeomorphic to $S^2$. 
By definition, {\it $\b$-plane} is a tangent null 2-plane at a point whose bivector 
is anti-self-dual, and {\it $\b$-surface} is a maximal embedded surface 
whose tangent plane is $\b$-plane at every point. 

LeBrun and Mason constructed the twistor correspondence for 
self-dual Zollfrei metrics. The statement is as follows; 
there is a natural one-to-one correspondence between 
\begin{itemize}
 \item equivalence classes of self-dual Zollfrei conformal structures 
   on $S^2\times S^2$, 
 \item equivalence classes of totally real embeddings 
    $\iota: \RP^3 \rightarrow\CP^3$, 
\end{itemize}
in neighborhoods of the standard conformal structure $[g_0]$ 
and the standard embedding of $\RP^3$. 
They also proved that the Zollfrei condition 
is an open condition in the space of self-dual neutral metrics. 
It implies that the term `Zollfrei' is removable in the above statement. 
This correspondence is characterized by the following condition; 
there is a double fibration $S^2\times S^2\overset{\wp}\longleftarrow\hat{\cZ} 
\overset{\Psi}\longrightarrow\CP^3$ such that 
\begin{enumerate}
 \item $\Psi$ is a continuous surjection and $\wp$ is a complex disk bundle, 
   i.e. for each $x\in S^2\times S^2$, $\hat{\cZ}_x=\wp^{-1}(x)$ is 
   biholomorphic to the complex unit disk, 
 \item $\Psi_x : \hat{\cZ}_x \rightarrow\CP^3$ is holomorphic on the interior of 
   $\hat{\cZ}_x$ and $\Psi_x(\del \hat{\cZ}_x)\subset P$, 
   where $\Psi_x$ is the restriction of $\Psi$ on $\hat{\cZ}_x$ 
   and $P=\iota(\RP^3)$,
 \item the restriction of $\Psi$ on $\hat{\cZ}-\del\hat{\cZ}$ 
   is bijective onto $\CP^2-N$, 
 \item $\{\wp(\Psi^{-1}(y))\}_{y\in P}$ is equal to the set of 
   $\b$-surfaces on $S^2\times S^2$. 
\end{enumerate}
The conditions (2) and (3) say that 
$\{\Psi(\hat{\cZ}_x)\}_{x\in S^2\times S^2}$ 
is a family of holomorphic disks on $\CP^3$ with boundaries on $P$ which foliate 
$\CP^3-P$. 

The twistor correspondence for the standard metric $g_0$ on 
$M=S^2\times S^2$ is explained in Lemma 8.1 of \cite{bib:LM05}. 
For the later convenience, we give an alternative description of the double fibration 
for $g_0$, and we describe the situation by using local coordinates. 
The twistor correspondence is described by the following diagrams 
which are explained below: 

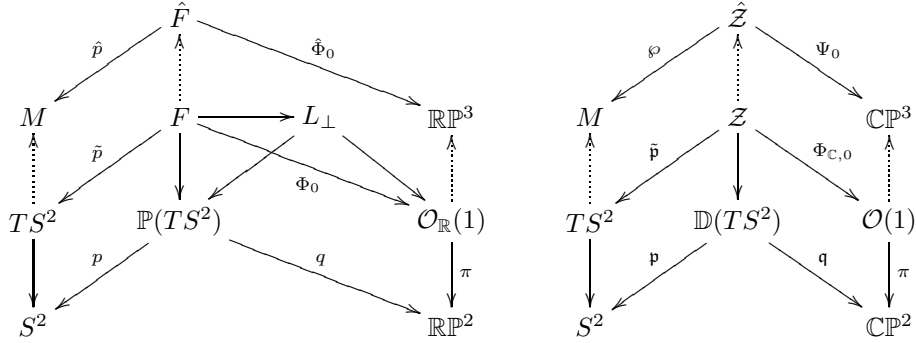
\begin{equation}\label{the_diagram}
 \xymatrix{
  & \hat{F} \ar[dl]_{\hat{p}} \ar[drr]^{\hat{\Phi}_0} 
   & & & & \hat{\cZ} \ar[dl]_{\wp} \ar[dr]^{\Psi_0} & \\
  M & F \ar@{.>}[u] \ar[dl]_{\tilde{p}} \ar[d] \ar[r] \ar[drr]_{\Phi_0} 
   & L_{\perp} \ar[dl] \ar[dr] & \RP^3 & M & 
   \cZ \ar@{.>}[u] \ar[dl]_{\tilde{\Gp}} \ar[d] \ar[dr]^{\Phi_{\C,0}} & \CP^3 \\
  TS^2 \ar@{.>}[u] \ar[d] & \bP(TS^2) \ar[dl]_p \ar[drr]^q & 
  & \cO_\R(1) \ar@{.>}[u] \ar[d]^\pi 
   & TS^2 \ar@{.>}[u] \ar[d] & {\mathbb D}(TS^2) 
     \ar[dl]_\Gp \ar[dr]^\Gq & \cO(1) \ar@{.>}[u] \ar[d]^\pi \\
  S^2 & & & \RP^2 & S^2 & & \CP^2 \\ }
\end{equation}
Here $p,q,\Gp$ and $\Gq$ are the same maps as in 
the diagrams (\ref{diag:double_fibration}), 
and the upward arrows are inclusions. 

It is convenient to use an identification of $S^2\times S^2$ with 
$\widetilde{Gr}_2(\R^4)$, where $\widetilde{Gr}_2(\R^4)$ is the Grassmannian 
consisting of oriented 2-planes in $\R^2$. 
Each element of $\widetilde{Gr}_2(\R^4)$ is represented by a $4\times 2$ matrix 
up to the right action of the group $GL_+(2,\R)$ consisting of $2\times 2$ 
matrices with positive determinant. We write $[\![a,b]\!]$ for the class 
represented by a $4\times 2$ matrix $(a,b)$ with rank two. 

Putting $y_0=[0:0:0:1]\in\RP^3$, 
$\RP^3-\{y_0\}$ has a line bundle structure over $\RP^2$ defined by 
$[z_1:z_2:z_3:z_4]\mapsto [z_1:z_2:z_3]$, 
and we denote this line bundle by $\cO_\R(1)$. 
Using the Euclidean metric on $\R^3$, $\cO_\R(1)$ is identified with the 
tautological bundle $\cL=\{([\z],v)\in\RP^2\times\R^3 : v\propto \z\}$ by 
$$ \cO_\R(1) \ni [t_1:t_2:t_3:\l] \longleftrightarrow 
  \([t_1:t_2:t_3], \l(t_1,t_2,t_3)\) \in \cL$$
where $t_1^2+t_2^2+t_3^2=1$. 

Let $F$ be the fiber product of the tangent bundle $TS^2\rightarrow S^2$ and 
$p : \bP(TS^2)\rightarrow S^2$, and let $\cO_\R(1)\rightarrow\RP^2$ be 
the tautological bundle. 
Let $F=L_\parallel\oplus L_\perp$ be the orthogonal decomposition 
over $\bP(TS^2)$ where $L_\parallel=\{(t,w,[v])\in F: w\propto v\}$ and
$L_\perp=\{(t,w,[v])\in F: w\perp v\}$. 
We define $\Phi_0 :F\rightarrow\cO_\R(1)$ to be the composition of 
the orthogonal projection $F\rightarrow L_\perp$ with 
the map $L_\perp \rightarrow \cL\simeq\cO_\R(1)$ 
given by $(t,w,[v])\mapsto ([t\times v],w)$. 

The embedding $TS^2\rightarrow M=\widetilde{Gr}_2(\R^4)$ is given by 
\begin{equation}\label{embedding_of_TS^2}
  TS^2 \ni (t,v) \longmapsto  
  \[\!\!\!\[ \begin{array}{cc} t_1 & -v_1 \\ t_2 & -v_2 \\ t_3 & -v_3 \\
  0 & 1 \end{array} \]\!\!\!\] \in \widetilde{Gr}_2(\R^4). 
\end{equation}
In this embedding, we have
$ TS^2 \cong\ \widetilde{Gr}_2(\R^4)-S^2_\infty$ where 
$ S^2_\infty $ consists of such a point given by
\begin{equation} \label{pt_in_S^2_infty}
 \[\!\!\!\[ \begin{array}{cc} t_1 & -v_1 \\ t_2 & -v_2 \\ t_3 & -v_3 \\
  0 & 0 \end{array} \]\!\!\!\]. 
\end{equation}

$TS^2$ is trivial over $D_\pm$ and $W$, and we define the coordinates 
$(x_1^\pm,x_2^\pm,x_3^\pm,x_4^\pm)\in\R^4\cong TS^2|_{D_\pm}$ and 
$(\a,\b,\e_1,\e_2)\in S^1\times[-\frac{\pi}{2},\frac{\pi}{2}]\times\R^2 
 \cong TS^2|W$ so as to fit  
\begin{equation}\label{coordinate_tangent}
 \[\!\!\!\[ \begin{array}{cc} \pm x_1^\pm & -x_4^\pm \\
 \pm x_2^\pm & x_3^\pm \\ \pm 1 & 0 \\ 0 & 1 
 \end{array} \]\!\!\!\] \quad \text{and} \quad
  \[\!\!\!\[ \begin{array}{cc} \cos\a & \e_2\sin\a \\
  \sin\a & -\e_2\cos\a \\ \tan\b & \e_1 \\
  0 & 1 \end{array} \]\!\!\!\], 
\end{equation}
in the manner of (\ref{embedding_of_TS^2}). 
Then the coordinate change is given by 
\begin{equation} \label{coordinate_change.II}
 \begin{cases}
   x_3^\pm = \e_1\sin\a\cot\b +\e_2\cos\a \\
   x_4^\pm = -\e_1\cos\a\cot\b + \e_2\sin\a. 
 \end{cases}
\end{equation}

The coordinates on $F$ is given by 
$(x_1^\pm,x_2^\pm,x_3^\pm,x_4^\pm,\z^\pm)$ and 
$(\a,\b,\e_1,\e_2,\xi)$, and the map $\Phi_0$ is described by 
\begin{equation}\label{map_Phi}
 \(x_1^\pm,x_2^\pm,x_3^\pm,x_4^\pm,\z^\pm\) \longmapsto 
   \[ 1 : \z^\pm : -x_1^\pm -x_2^\pm \z^\pm : 
        -x_3^\pm \z^\pm + x_4^\pm \], \end{equation}
$$ \(\a,\b,\e_1,\e_2,\x\) \longmapsto \[ \hasira  
 \x\cos\a\tan\b-\sin\a : -\x\sin\a\tan\b +\cos\a : \x : -\e_1\x+\e_2 \]. $$

Let $\cO(1)\cong \CP^3-\{y_0\}$ be the complex line bundle defined in the 
similar way as $\cO_\R(1)$, and $\cZ$ be the fiber product of 
$TS^2\rightarrow S^2$ and $\Gp : \cD(TS^2)\rightarrow S^2$. 
Similar to the Zoll case, we obtain the double fibration 
$TS^2\leftarrow\cZ\rightarrow\cO(1)$ 
by extending real parameters $\z^\pm$ and $\x$ to the complex parameters. 

The double fibration 
$M=\widetilde{Gr}_2(\R^4)\leftarrow\hat{\cZ}\rightarrow\CP^3$ 
is obtained as follows. 
For each $x\in M$, we define a holomorphic disk $D_x$ in 
$\CP^3$ with boundary on $\RP^3$ by the following; 
$D_x=\Phi_{\C,0}(\Gp^{-1}(x))$ if $x\in TS^2$, and 
$$ D_x= \{ [z_1:z_2:z_3:u] : \Imag u \geq 0 \} \cup \{ y_0 \}$$ 
if $x\in S^2_\infty$ of the form (\ref{pt_in_S^2_infty}), 
where $(z_1,z_2,z_3)=t\times v$.
Then $\{D_x\}_{x\in M}$ is a family of holomorphic disks in $(\CP^3,\RP^3)$ 
which foliates $\CP^3-\RP^3$. Let
$$ \hat{\cZ}=\{ (x,y)\in M\times\CP^3 : y\in D_x\}, $$
then $\hat{\cZ}$ has a natural smooth structure such that the maps 
$\wp : \hat{\cZ}\rightarrow M$ and $\Psi : \hat{\cZ}\rightarrow\CP^3$ are smooth. 
In this way, we obtain the double fibration $M\leftarrow\hat{\cZ}\rightarrow\CP^3$. 
The real version $M\leftarrow\hat{F}\rightarrow\RP^3$ is obtained by 
restricting $\hat{\cZ}$ to the boundary. \\

\noindent
{\it Remark :} By the identification 
$\widetilde{Gr}_2(\R^4)\cong S^2\times S^2$, 
the standard neutral metric $g_0$ is described in the coordinates of 
$TD_\pm$ by
$$ \frac{1}{1+ \| x\|^2 + \D^2}
   \(dx_1^\pm dx_3^\pm + dx_2^\pm dx_4^\pm \) 
   \quad \text{where} \ 
   \begin{cases} 
   \ \| x\|^2= \sum_{i=1}^4 (x_i^\pm)^2, \\
   \ \D=x_1^\pm x_3^\pm + x_2^\pm x_4^\pm. 
   \end{cases} $$

\noindent
{\it Remark :} The point $y=[z_1:z_2:z_3:z_4]\in\RP^3$ corresponds to the 
embedded two sphere $\{ [\![ a, b ]\!] : ya=yb=0 \} \subset 
\widetilde{Gr}_2(\R^4)$ which is a $\b$-surface for $g_0$.  
On the other hand, the holomorphic disk corresponding to 
$[\![ a,b ]\!]\in\widetilde{Gr}_2(\R^4)$ is the closure of one of the 
two connected components of $\{y\in\CP^3-\RP^3 : ya=yb=0 \}$. 

\section{Definition of the singularity} \label{Section:Def_sing}

In this section, we introduce a certain singularity of metrics and 
twistor spaces, 
and state a conjecture for a singular version of the twistor correspondence. 
The main purpose in this article is to construct 
explicit examples such that the conjecture holds restricting to these examples. 
That is explained in the following sections. 

First, we study about the twistor spaces for the self-dual Zolfrei metrics. 
In the non-singular case, as explained in Section \ref{Section:Standard}, 
the twistor space is given by a pair 
$(\CP^3,P)$ where $P$ is the image of a totally real embedding 
$\iota : \RP^3\rightarrow \CP^3$. 
In \cite{bib:LM05}, 
C.~LeBrun and L.~J.~Mason proved that if the embedding $\iota$ 
is close to the standard one, then there is a unique smooth family of 
holomorphic disks in $\CP^3$ with boundaries on $P$ which satisfy  
(1) the relative homology class of each disk generates $H_2(\CP^3,P;\Z)\cong\Z$, 
and (2) interiors of these disks smoothly foliate $\CP^3-P$. 

We formulate the singular version of the above conditions. 
Let $\iota:\RP^3\rightarrow\CP^3$ be the continuous injection such that 
the restriction $\RP^3-\{y_0\}\rightarrow \CP^3$ is a totally real embedding, 
and let $P=\iota(\RP^3)$. 
We assume $\iota$ is homotopic to the standard embedding 
$\RP^3\subset\CP^3$, then we have $H_2(\CP^3,P; \Z)\cong\Z$. 
\begin{Def}
 The pair $(\CP^3,P)$ is said to satisfy the condition $(\sharp)$ if there is 
 a unique family of holomorphic disks in $\CP^3$ with boundaries on $P$ 
 which satisfy 
 \begin{list}{($\sharp$\,\arabic{mynum})}%
 {\usecounter{mynum} \itemsep 0in  \leftmargin .4in}
 \item the relative homology class of each disk generates $H_2(\CP^3,P;\Z)$,
 \item interiors of these disks foliate $\CP^3-P$, and the disks which do 
     not contain $y_0$ forms a smooth family.  
\end{list} 
\end{Def}
\noindent
In Sections \ref{Section:twistor} and \ref{Section:holo_disks}, 
we construct examples of the pair $(\CP^3,P)$ 
which is given by such singular embeddings 
and which satisfy the condition $(\sharp)$. 

Next, we define the ``singular self-dual Zollfrei metric". 
Let $M$ be a four dimensional smooth manifold, $C$ be 
a two dimensional closed submanifold of $M$, 
and $g$ be a neutral metric on $M-C$. 
\begin{Def} \label{Def:sing_b-surf}
 $g$ is called a {\it singular neutral metric} with {\it singular $\b$-surface} $C$ 
 if, for all $x\in C$, there is an open neighborhood 
 $x\in U\subset M$ with the coordinate $u=(u_1,u_2,u_3,u_4)$ on $U$, 
 which satisfies following two conditions: 
 \begin{enumerate}
  \item $C\cap U=\{u\in U : u_2=u_3=0 \}$. 
  \item Let $ j :(u_1,u_4,r,\phi) \mapsto (u_1,u_2,u_3,u_4)$ be the 
   cylindrical coordinate given by $u_2=r\cos\phi, u_3=r\sin\phi$. 
   Then, for some smooth function $h$ on $U\setminus C$,
   we can write $hg=g_{\rm st}+g_1+o(r)$ such that 
   $g_{\rm st}=2(du_1du_3+du_2du_4)$ is the standard neutral metric, 
   $g_1$ is a symmetric tensors satisfying
   $j^*g_1=\rho(u_1,u_4,\phi)r^2d\phi^2$ for some function $\rho$, 
   and $o(r)$ is the error term satisfying $\lim_{r\rightarrow 0}o(r)=0$. 
 \end{enumerate}
\end{Def}
\noindent
We can also define the {\it singular neutral conformal structure} on $M$ 
with singular $\b$-surface $C$. 

\vspace{1ex}
\noindent
{\it Remark :} If we take the limit $r\rightarrow 0$ for fixed $(u_1,u_4,\phi)$, 
then $g_0+g_1+o(r) \rightarrow g_0+
\rho\cdot(-\sin\phi du_2+\cos\phi du_3)^2$. 
This limit defines a neutral metric on $T_xM$ where $x=\lim u\in C$. 
This metric depends on $\phi$, but $\{u_2=u_3=0\}$ always defines a $\b$-plane. 
This is the reason why we call $C$ ``singular $\b$-surface". 

\begin{Def} \label{Def:sing_SD_Zf}
 (i) A neutral metric $g$ on $M-C$ is 
 called a {\it singular self-dual neutral metric} on $M$ when 
 \begin{enumerate}
  \item $g$ is self-dual on $M-C$, and 
  \item $g$ is a singular neutral metric with singular $\b$-surface $C$. 
 \end{enumerate}
  (ii) $g$ is called a {\it singular Zollfrei metric} when 
  every null geodesic in $M-C$ satisfies either of the followings; 
  it is closed or the ends of its closure in $M$ are the points in $C$ 
  which are not necessary distinct. 
\end{Def}

In the non-singular case, every $\b$-surface is either $S^2$ or $\RP^2$ 
if the metric $g$ is self-dual Zollfrei (Theorem 5.14 of \cite{bib:LM05}). 
In our case, $\b$-surface is defined only on $M-C$. However 
we will see later that, in our examples, the closure of each $\b$-surface in $M$ 
is homeomorphic to $S^2$, with extra two points
(Proposition \ref{extension_of_b-surface} and \ref{Prop:conjecture_holds}). 
Motivated by these examples we state the following conjectures. 
\begin{Conj} \label{Conj:b-surface}
 Let $g$ be a singular self-dual Zollfrei metric on $M$ with a singular 
 $\b$-surface $C$. 
 Let $S$ be any $\b$-surface in $M-C$, and $\bar{S}$ be its closure in $M$. 
 Then $\bar{S}-S$ is a finite subset of $C$ 
 and $\bar{S}$ is a topological manifold. 
\end{Conj}
\noindent
Let $g$ and $S$ be as in the above conjecture. 
We simply call $\bar{S}$ as a $\b$-surface for $g$. 
\begin{Conj} \label{Conj:twistor_corr}
 There is a natural one-to-one correspondence between 
 \begin{itemize}
  \item equivalence classes of singular self-dual Zollfrei conformal structures on 
    $S^2\times S^2$ with a single singular $\b$-surface $C\simeq S^2$, 
  \item equivalence classes of the pairs $(\CP^3,P)$ which satisfy the condition 
  $(\sharp)$, where $P$ is the image of a totally real embedding 
  $\RP^3\rightarrow\CP^3$ which has one singular point.  
 \end{itemize}
\end{Conj}
\noindent 
This correspondence should be characterized by the similar condition as 
LeBrun and Mason's theorems, and the explicit formulation for our examples is given in 
Theorem \ref{Thm:twistor_correspondence}. 
We formulated the conjecture for the simplest possible singularity 
of the twistor space 
because otherwise more complicated phenomenon would occur such as 
intersection of two singular $\b$-surfaces. 

\section{Construction of singular metrics} \label{Section:sing_met}
In this section, we construct the examples of singular self-dual Zollfrei metric, 
by patching a pair of Petean's indefinite self-dual metrics. 
We will see later that each example corresponds to an element of 
the set $\cS(\R^2)^{\text{sym}}$ consisting of functions that are 
rapidly decreasing, axisymmetric, and real valued. 
Recall that a smooth function $f(x)$ on $\R^n$ is called {\it rapidly decreasing} 
if and only if for each polynomial $P$ and each integer $m\geq 0$ 
$$ \sup \left| \hasira |x|^m P(\del_1,\cdots,\del_n)f(x)\right| <\infty $$ 
where $|x|$ denote the norm of $x$. 
We write $\cS(\R^n)$ for the set consisting of rapidly decreasing 
real valued functions on $\R^n$. 
We call $f\in\cS(\R^2)$ {\it axisymmetric} if and only if $f$ is 
$SO(2)$-invariant. 

{\it Petean's metric} is an indefinite metric over $\R^4$ of the form
\begin{equation}\label{Petean_met}
 g=2(dx_1 dx_3+dx_2 dx_4) + f(x_1,x_2) \( dx_1^2+dx_2^2 \), 
\end{equation}
where $f$ is a smooth function.
Such metrics are first introduced by J.~Petean to construct examples 
of indefinite K\"{a}hler-Einstein metric on the complex tori or 
the primary Kodaira surfaces (\cite{bib:Petean}). 
For the metric (\ref{Petean_met}), 
we have an indefinite orthonormal frame $\{e_i\}_i$ given by 
$$ (e_1,e_2,e_3,e_4)=\(\del_1,\del_2,\del_3,\del_4\)
  \frac{1}{2\sqrt{D}}
  \( \begin{array}{rrrr} 1&-1&-1&1 \\ 1&1&-1&-1 \\ -b&b&a&-a \\ -b&-b&a&a 
  \end{array} \), $$
where $\del_i=\pd{x_i}$ and 
$$ D=f^2+4, \ \ 
   a=\frac{f+\sqrt{D}}{2}, \ \ b=\frac{f-\sqrt{D}}{2}. $$ 
Let $\wedge^2_- T\R^4$ be the bundle of anti-self-dual bivectors. 
If we trivialize $\wedge^2_- T\R^4$ by the frame
$$ e_1\wedge e_2 - e_3\wedge e_4, \ \ 
   e_1\wedge e_3 - e_2\wedge e_4, \ \ 
   e_1\wedge e_4 + e_2\wedge e_3, $$
then we can check that the induced Levi-Civita connection is trivial.

\noindent
{\it Remark :} Petean's metric is a special case of the {\it Walker metric} 
(cf.\cite{bib:Matsushita04},\cite{bib:Matsushita05}). 
And the orthonormal frame defined above looks similar to the one explained 
in \cite{bib:Matsushita04}. 

Following the arguments in \cite{bib:LM05}, the vectors 
\begin{equation}\label{m_1m_2}
 \begin{aligned}
   \Gm_1&=e_1-\sin 2\s\ e_3 + \cos 2\s\ e_4 \\
   \Gm_2&=e_1+\cos 2\s\ e_3 + \sin 2\s\ e_4 
 \end{aligned}\end{equation}
span $\b$-plane at every point in $\R^4$ for each $\s\in\RP^1=\R/\pi\Z$. 
If we put
\begin{equation}\label{n_1n_2}
 \begin{aligned}
   \Gn_1&=\cos\s\ \del_1 + \sin\s\ \del_2
          -\frac{f}{2}\( \cos\s\ \del_3 + \sin\s\ \del_4 \), \\
   \Gn_2&=\sin\s\ \del_3 - \cos\s\ \del_4,  
 \end{aligned}\end{equation}
then the distribution $\langle\Gm_1,\Gm_2\rangle$ is equal to 
$\langle\Gn_1,\Gn_2\rangle$. 

\begin{Prop}\label{Prop:b-surface}
 Let $g$ be a Petean's metric of the form (\ref{Petean_met}), 
 then every $\b$-surface on $(\R^4,g)$ 
 is given by the solutions of 
 \begin{equation}\label{eq:b-surface}
  \begin{cases}
    -\sin\s\ x_1+\cos\s\ x_2= c_1, \\
    \cos\s\ x_3+\sin\s\ x_4+ \vp(x_1,x_2,\s)=c_2, 
   \end{cases}
 \end{equation}
 for some real constants $\s, c_1$ and $c_2$, 
 where
 \begin{equation}\label{eq:vp}
 \vp(x_1,x_2,\s)= \frac{1}{2} \int_0^\l 
  f(\cos\s\ t-\sin\s\ \mu,\ \sin\s\ t+ \cos\s\ \mu) dt, 
\end{equation}
$$  \( \begin{aligned} \l \\ \mu \end{aligned} \) = 
    \( \begin{array}{rr} \cos\s & \sin\s \\ -\sin\s & \cos\s 
	  \end{array} \) 
	\( \begin{aligned} x_1 \\ x_2 \end{aligned} \). $$
\end{Prop}
\begin{proof}
Notice that $\vp$ satisfies 
\begin{equation}
  \cos\s\ \PD{\vp}{x_1} + \sin\s\ \PD{\vp}{x_2} = \frac{f}{2}.
\end{equation}
So the left hand sides of (\ref{eq:b-surface}) are both annihilated by 
${\mathfrak n}_1$ and ${\mathfrak n}_2$. Hence these are constant 
along some $\b$-surface. 
\end{proof}
Because $(\s,c_1,c_2)$ and $(\s+\pi,-c_1,-c_2)$ correspond to the same 
$\b$-surface, we can assume $\s\in[0,\pi)$. 

\vspace{1ex}
\noindent
{\it Remark :} In \cite{bib:BDM}, 
D.~E.~Blair et.al. showed that the `hyperbolic twistor space' over $\R^4$ 
equipped with a Petean's metric is holomorphically trivial. 
They proved this fact by constructing an explicit complex coordinate 
of the twistor space. 
This construction actually works in the case of the `reflector space' (\cite{bib:JR})
or, in the other literature, the `product twistor space' (\cite{bib:Blair}). 
In this way, we obtain essentially the same statement as 
Proposition \ref{Prop:b-surface}. 

\vspace{1ex}
\noindent
{\it Remark :} One can show a Petean's metric in the form (\ref{Petean_met}) is 
flat if and only if $f$ is harmonic (i.e. $f_{x_1x_1}+f_{x_2x_2}=0$).  
Hence, for example, if we assume $f$ is rapidly decreasing, 
then the metric is flat if and only if $f\equiv 0$. 
(cf.\cite{bib:Matsushita05},\cite{bib:Petean})

\begin{Def}
Let $g$ be a Petean's metric of the form (\ref{Petean_met}), then $g$ is called 
{\it rapidly decreasing} (respectively {\it compact supported}, or 
{\it axisymmetric}) if and only if $f$ is so. 
On the other hand, the {\it dual of $g$} is another Petean's metric 
of the form
$$ g^\vee =2(dx_1 dx_3+dx_2 dx_4) - f(x_1,x_2) \( dx_1^2+dx_2^2 \). $$
\end{Def}

\vspace{1ex}
\noindent
{\it Remark :} In the space of Petean's metrics of the form (\ref{Petean_met}), 
the flat metric, i.e. the case of $f\equiv 0$, is characterized by 
an $(S^1\times S^1)$-invariance, 
where $(\t_1,\t_2)\in S^1\times S^1$ acts on $\R^4$ by 
$$ \begin{pmatrix} x_1 & -x_4 \\ x_2 & x_3 \end{pmatrix}
    \longmapsto R(\t_1) 
	\begin{pmatrix} x_1 & -x_4 \\ x_2 & x_3 \end{pmatrix} R(\t_2)^{-1}, 
	\quad \( R(\t)=\begin{pmatrix} \cos\t & -\sin\t \\ \sin\t & \cos\t 
	       \end{pmatrix} \) .$$
In the same way, the axisymmetric metrics are characterized by 
the $(S^1\times\{1\})$-invariance. 
On the other hand, the `dual' defines a $\Z_2$-action on the space of Petean's 
metrics, then the flat metric is also characterized by the $\Z_2$-invariance. 

\vspace{1ex}
We use the coordinates introduced in Section \ref{Section:Standard} from now on. 
\begin{Prop} \label{extension_of_b-surface}
  Let $g_\pm$ be compact supported Petean's metrics on $TD_\pm$ given by 
  \begin{equation} \label{metrics_on_TD_pm} 
    g_\pm=2(dx_1^\pm dx_3^\pm+dx_2^\pm dx_4^\pm) +
      f_\pm(x_1^\pm,x_2^\pm) \( (dx_1^\pm)^2+(dx_2^\pm)^2 \). 
  \end{equation}
  Then the conformal structures of these metrics extend 
  to a self-dual indefinite conformal structure on $TS^2$. 
  Moreover, if $g_\pm$ are both axisymmetric and dual each other, 
  then every $\b$-surface on $TS^2$ is an embedded $S^1\times\R$ 
  whose closure in $\widetilde{Gr}_2(\R^4)$ is homeomorphic to $S^2$, 
  with extra two points. 
\end{Prop}
\begin{proof}
If $g_\pm$ are flat i.e. $f_\pm\equiv 0$, it is obvious that these conformal 
structures extend to the standard conformal structure on 
$\widetilde{Gr}_2(\R^4)$ by the remark in Section \ref{Section:Standard}. 
In the general case, 
because we assumed $f_\pm$ are compact supported, 
these metrics are flat in some neighborhood of 
$TS|_{W_0}$ where $W_0=\{(\a,\b)\in W: \b=0 \}.$
So these conformal structures extend to a self-dual conformal structure on $TS^2$. 

Now we suppose the axisymmetricity and the duality, 
and denote simply $f_+=-f_-=f$. 
Because $f$ is compact supported, we can take $R>0$ 
such that $f(x_1^+,x_2^+)=0$ when 
$(x_1^+)^2+(x_2^+)^2>R^2$. 
By the Proposition \ref{Prop:b-surface}, each $\b$-surface on $TD_\pm$ is given by 
 \begin{equation} \label{eq:beta-surfaces}
  \begin{cases}
    -\sin\s^\pm\ x_1^\pm+\cos\s^\pm\ x_2^\pm= c_1^\pm, \\
    \cos\s^\pm\ x_3^\pm+\sin\s^\pm\ x_4^\pm \pm 
	   \vp(x_1^\pm,x_2^\pm,\s^\pm)=c_2^\pm, 
   \end{cases}
 \end{equation}
for some $(\s^\pm,c_1^\pm,c_2^\pm)$, where $\vp$ is 
defined by (\ref{eq:vp}). 
We observe that when $(\s^+,c_1^+,c_2^+)=(\s^-,c_1^-,c_2^-)$, 
corresponding $\b$-surfaces are nicely extended in $TS^2$. 
So we drop the signs on $\s,c_1,c_2$ for the simplicity. 
Changing the coordinates by (\ref{coordinate_change}) and 
(\ref{coordinate_change.II}), the $\b$-surfaces are described in the regions 
$\{(\a,\b,\e_1,\e_2)\in TW :0< |\tan\b| <R^{-1}\}$ by
\begin{equation}\label{coordinate_changed_eq}
 \begin{cases} 
  \sin(\a-\s)=c_1 \tan\b, \\
  -c_1\e_1-\cos(\a-\s)\e_2 +\psi(\a,\b,\s) =c_2, 
  \end{cases}
\end{equation}
where $\psi$ is defined as the following. Let
$$ \hat{f}^{(\pm)}(\s,\mu)=
 \int_0^{\pm\infty} f(\cos\s\ t-\sin\s\ \mu,\ \sin\s\ t+ \cos\s\ \mu) dt $$
be the `half Radon transform' of $f$. Then we define $\psi$,  
when $0<\tan\b<R^{-1}$, 
$$ \psi(\a,\b,\s)= \begin{cases}
   \frac{1}{2} \hat{f}^{(+)}(\s,\sin(\a-\s)\cot\b) 
     & \qquad \text{if} \quad \cos(\a-\s)>0 \\
   \frac{1}{2} \hat{f}^{(-)}(\s,\sin(\a-\s)\cot\b)
     & \qquad \text{if} \quad \cos(\a-\s)<0 \end{cases} $$
and when $-R^{-1}<\tan\b<0$, 
$$ \psi(\a,\b,\s)= \begin{cases}
   -\frac{1}{2} \hat{f}^{(-)}(\s,\sin(\a-\s)\cot\b)
     & \qquad \text{if} \quad \cos(\a-\s)>0. \\
	 -\frac{1}{2} \hat{f}^{(+)}(\s,\sin(\a-\s)\cot\b) 
     & \qquad \text{if} \quad \cos(\a-\s)<0 \end{cases} $$
Let $\hat{f}$ be the Radon transform which is explained in 
Section \ref{Appendix}, then we have 
$\hat{f}^{(+)}=-\hat{f}^{(-)}=\frac{1}{2}\hat{f}$ from the axisymmetricity. 
Moreover $\hat{f}(\mu)=\hat{f}(\s,\mu)$ is even for $\mu$ and 
does not depend on $\s$. 
So the equation (\ref{coordinate_changed_eq}) extends to $|\tan\b|<R^{-1}$ 
by exchanging $\psi(\a,\b,\s)$ with
\begin{equation}\label{tilde_psi} 
  \tilde\psi(\a,\s;c_1)=\begin{cases}
   \frac{1}{4} \hat{f}(c_1) 
     & \qquad \text{if} \quad \cos(\a-\s)>0 \\
   -\frac{1}{4} \hat{f}(c_1)
     & \qquad \text{if} \quad \cos(\a-\s)<0. \end{cases}
\end{equation}
Notice that $\cos(\a-\s)\neq 0$ near $\{\b =0\}$ because 
$\sin(\a-\s)$ is close to zero by the first equation of 
(\ref{coordinate_changed_eq}). 
Hence we have proved that a suitable pair of $\b$-surfaces on $TD_\pm$  
extends to a $\b$-surface on $TS^2$.
This $\b$-surface has a fiber bundle structure with respect to the projection 
$TS^2\rightarrow S^2$. In fact the base space is a big circle on $S^2$ 
and the fibers are orientable linear lines in some tangent space of $S^2$. 
So it is isomorphic to a cylinder $S^1\times \R$. 
If we take any point on this $\b$-surface, and let it move 
away along the fiber, then the point goes to
\begin{equation} \label{antipodal_points}
  \[\!\!\!\[ \begin{array}{cc} 
  \mp c_1\sin\s & \cos\s \\ \pm c_1\cos \s & \sin\s \\ \pm 1 & 0 \\
  0 & 0 \end{array} \]\!\!\!\] \in S^2_\infty, \end{equation}
so the proof is completed.
\end{proof}

The two points (\ref{antipodal_points}) are actually {\it antipodal}, 
i.e. they are exchanged by the natural involution of $S^2_\infty$ 
which is defined as a deck transformation
for the orientation forgetting map 
$\widetilde{Gr}_2(\R^4) \rightarrow Gr_2(\R^4)$. 

\begin{Prop} \label{Prop:singularity_of_metric}
 Let $g_\pm$ be Petean's metrics on $TD_\pm$ in the form 
 (\ref{metrics_on_TD_pm}), and we assume that $g_\pm$ are 
 compact supported, axisymmetric and dual each other. 
 Then the induced self-dual conformal structure on 
 $TS^2=\widetilde{Gr}_2(\R^4)-S^2_\infty$ defines a singular self-dual conformal 
 structure on $\widetilde{Gr}_2(\R^4)$ in the sense of Definition \ref{Def:sing_b-surf}. 
\end{Prop} 
\begin{proof}
 We check that $S^2_\infty$ is a singular $\b$-surface. 
 As an example, we study on the coordinate neighborhood $U$ given by 
 $$ \R^4 \ni (u_1,u_2,u_3,u_4) \longmapsto 
  \[\!\!\!\[ \begin{array}{cc} u_1 & -u_4 \\ 1&0  \\ 0&1 \\ u_2 & u_3
 \end{array} \]\!\!\!\] \in \widetilde{Gr}_2(\R^4). $$
 In this coordinate, $S^2_\infty$ is given by $\{u_2=u_3=0\}$. 
 We can check that $U\cap TD_+$ is given by 
 $\{x\in TD_+: x_3^+<0\}=\{u\in U: u_2<0\}$, 
 and by a direct calculation, 
 \begin{Eqnarray*} 
  g_+&=&\frac{1}{u_2^2} \left[ 2(du_1du_3+du_2du_4) \hasira \right. \\
     && \quad  \left. \hasira + r^2 f(u_1\tan\phi-u_4,-\tan\phi)
	 \left\{(1+u_1^2)d\phi^2+(\sin\phi du_1+\cos\phi du_4)^2 \right\} \right]
 \end{Eqnarray*}
 where $u_2=r\cos\phi, u_3=r\sin\phi$.  
 If we put $\rho=f\cdot (1+u_1^2)$ and so on, the condition (2) of 
 Definition \ref{Def:sing_b-surf} is satisfied. 
 
 In this way, checking several coordinate neighborhoods, 
 we can show that $S^2_\infty$ is a singular $\b$-surface for the metric. 
\end{proof}

\noindent
{\it Remark :} The $(S^1\times S^1)$-action on $TD_\pm$ extends to  
$\widetilde{Gr}_2(\R^4)$ by
$$ [\![a,b]\!] \longmapsto [\![ \tilde{R}(\t_1,\t_2) (a,b) ]\!]
   \qquad \( \tilde{R}(\t_1,\t_2)= 
   \begin{pmatrix} R(\t_1) & O \\ O & R(\t_2) \end{pmatrix} \). $$
The standard metric on $\widetilde{Gr}_2(\R^4)$ is invariant under this action, 
and our singular metric is invariant under $(S^1\times \{1\})$-action. \\

Next we check the Zollfrei condition for our singular metric. 
In the non-singular case, self-dual metric $g$ is Zollfrei if every $\b$-surface 
is an embedded $S^2$ (Theorem 5.14. of \cite{bib:LM05}). 
Although we can not apply this theorem in our case, 
we can check the condition directly by writing down null geodesics explicitly. 

First, the following formulas are checked by a direct calculation: 
\begin{Lem} \label{Lem:nabla}
 Let $\nabla$ be the Levi-Civita connection of a Petean's metric 
 of the form (\ref{Petean_met}) and let $D=f^2+4$, then
  \begin{equation} \begin{aligned}
    & \nabla_{\del_3} =\nabla_{\del_4}=0, \\
	& \nabla_{\del_1}\del_3 =\frac{f\del_1f}{2D}\del_3, \quad
	 \nabla_{\del_1}\del_4 =\frac{f\del_1f}{2D}\del_4, \\
    & \nabla_{\del_1}\del_1 =\frac{f\del_1f}{2D}\del_1 
	      +\frac{\del_1f}{2} \del_3 -\frac{\del_2f}{2}\del_4. 
  \end{aligned} \end{equation}
\end{Lem}

\begin{Prop} 
 The singular self-dual metric on $\widetilde{Gr}_2(\R^4)$ in 
 Proposition \ref{Prop:singularity_of_metric} is singular Zollfrei. 
\end{Prop}
\begin{proof}
 Because every null geodesic is contained in some $\b$-surface, 
 the image of a null geodesic by the projection $TS^2\rightarrow S^2$ 
 is contained in some big circle of $S^2$. 
 We prove that this image is either one point or a whole circle. 
 If the image is one point, from (\ref{eq:beta-surfaces}) or
 (\ref{coordinate_changed_eq}), 
 the null geodesic must be a linear line in some tangent space of $S^2$. 
 This is actually a geodesic, 
 because, for example on $TD_+$, we have 
 $\nabla_{\del_3}=\nabla_{\del_4}=0$ by Lemma \ref{Lem:nabla}. 
 Notice that the end points of such a linear line are the antipodal points 
 on $S^2_\infty$. 

 For the other case, we first study about a null geodesic in $TD_+$ 
 whose image on $D_+$ is not one point. 
 Without any loss of generality, 
 we can assume this geodesic is contained in a $\b$-surface 
 with $\s^+=0$ in (\ref{eq:beta-surfaces}) because of the axisymmetricity. 
 Let $c(s)$ be a curve contained in this $\b$-surface. If the projected image 
 on $D_+$ is not one point, such a curve is always written in the form
 $$ x_1^+=s,\ \ x_2^+=c_1, \ \
   x_3^+=c_2-\vp (s,c_1,0), \ \
   x_4^+=\nu(s), $$
 at least for a small interval of parameter $s$, 
 where $\nu(s)$ is an unknown function. 
 Then the velocity vector is 
 $$ \dot{c}(s)=\del_1 - \frac{f}{2} \del_3 + \PD{\nu}{s} \del_4. $$
 Using Lemma \ref{Lem:nabla}, we have 
 $$ \tilde{\nabla}_{\frac{\del}{\del s}} \dot{c}(s) = \frac{f\del_1f}{2D} 
 \dot{c}(s) + \(\frac{d^2\nu}{ds^2}-\frac{\del_2f}{2}\) \del_4, $$
 where $\tilde{\nabla}$ is the covariant derivative along $c(s)$.  
 Hence $c(s)$ is an unparameterized geodesic if and only if 
 \begin{equation} \label{ODE_nu}
   \frac{d^2\nu}{ds^2}=\frac{1}{2}\del_2f(s,c_1). 
 \end{equation}
 Let $\nu_0(s)$ be the solution of (\ref{ODE_nu}) with $\nu_0(0)=\nu'_0(0)=0$, 
 then any solution of (\ref{ODE_nu}) is given by 
 $$\nu(s)=\nu_0(s)+q_1s+q_2, $$
 for some constants $q_1$ and $q_2$. 
 Every null geodesic on $TD_+$ whose image to $D_+$ is not one point
 is given by rotating the above $c(s)$. 
 Remark that $\nu_0$ is an even function, 
 moreover, because $f$ is compact supported, $\nu_0$ is a degree one polynomial
 for $|s|\gg 0$. So we have
 \begin{equation} \label{nu_0_at_large}
   \nu_0(s)= A_1|s|+A_2 \qquad \(|s|>R\),
 \end{equation}
 where $A_1$ and $A_2$ are constants, and
 $R$ is a large constant such as $f(x)\equiv 0$ for $|x|>R$. 

 Now notice the following geodesics on $TD_\pm$ parameterized by 
 $s_\pm$ respectively: 
 $$ x_1^\pm=s_\pm,\ \ x_2^\pm=c_1, \ \
   x_3^\pm=c_2\mp\vp (s_\pm,c_1,0), \ \
   x_4^\pm=\pm\nu_0(s_\pm)+q_1^\pm s_\pm + q_2^\pm. $$
 We prove that these geodesics extend smoothly in $TS^2$ 
 for suitable $q_1^\pm$ and $q_2^\pm$, 
 by changing the parameter 
 $$ u_+= \begin{cases} -\frac{1}{s_+} & s_+>0 \\
        -\frac{1}{s_-} & s_-<0, \end{cases} \qquad 
	u_-= \begin{cases} -\frac{1}{s_-} & s_->0 \\
    	-\frac{1}{s_+} & s_+<0. \end{cases} $$
 Changing the coordinates by (\ref{coordinate_change}) and 
 (\ref{coordinate_change.II}), these geodesics are written 
 on $u_+<0$
 $$ \begin{aligned}
  \a&=\tan^{-1}(-c_1u_+) \quad \(-\frac{\pi}{2}\leq\a\leq\frac{\pi}{2}\), 
     \quad \b=\tan^{-1}\(-u_+(1+c_1^2u_+^2)^{-\frac{1}{2}}\), \\
  \e_1&=\frac{1}{1+c_1^2u_+^2} 
       \( c_1u_+^2 \Theta + \Xi \), 
   \qquad 
     \e_2=\frac{1}{(1+c_1^2u_+^2)^\frac{1}{2}}
       \( \Theta -c_1\Xi \), 
  \end{aligned} $$
 where
 $$ \Theta=c_2-\vp\(-u_+^{-1},c_1,0\), \quad
      \Xi=u_+\nu_0\(-u_+^{-1}\)-q_1^++q_2^+u_+. $$
 Similarly on $u_+>0$, $(\a,\b,\e_1,\e_2)$ are given by the same equation
 with exchanging 
 $$ \Theta=c_2+\vp\(-u_+^{-1},c_1,0\), \quad
      \Xi=-u_+\nu_0\(-u_+^{-1}\)-q_1^-+q_2^-u_+. $$
 By the similar argument in the proof of Proposition \ref{extension_of_b-surface},
 $\Theta$ extends smoothly to $u_+=0$, and from
 (\ref{nu_0_at_large}), $\Xi$ extends smoothly to $u_+=0$ iff 
 $q_1^-= q_1^+$ and $q_2^-=q_1^++2A_1$. 
 By the similar argument for $u_-$, we obtain the same conditions. 
 So these geodesics are nicely extended to a closed curve if 
 $q_1^-= q_1^+$ and $q_2^-=q_1^++2A_1$. 
 
 The rest possibility is that the projected image of a null geodesic 
 is contained in the equator $W_0$. 
 However these null geodesics are of course closed, 
 because the conformal structure is standard around $TS^2|_{W_0}$. 
\end{proof}

So far, we assumed that the Petean's metrics are all compact supported, 
but this assumption is weakened to be `rapidly decreasing'. 
Actually the argument of  the extension to $TS^2|_{W_0}$ works well 
essentially in the same manner, and only one thing that we have to check is 
the smoothness at infinity, 
which is almost obvious from the rapidly decreasing condition. 
Getting together, we have 
\begin{Thm} \label{Thm:Sing_met}
Let $g_+$ and $g_-$ be rapidly decreasing Petean's metrics on $\R^4$, 
Then the disjoint union $(\R^4,[g_+]) \amalg (\R^4,[g_-])$
naturally extends to a self-dual indefinite conformal structure on 
$\widetilde{Gr}_2(\R^4)-S^2_\infty$. 
Moreover, if $g_\pm$ are both axisymmetric and dual each other, 
then this conformal structure defines a singular self-dual 
Zollfrei conformal structure on $\widetilde{Gr}_2(\R^4)$ 
with singular $\b$-surface $S^2_\infty$. 
\end{Thm}

\noindent
{\it Remark :}
These are the required examples. 
Notice that each singular self-dual Zollfrei conformal structure corresponds to 
a function $f\in \cS(\R^2)^{\text{sym}}$ by 
using the Petean's metric  $g_+$ which corresponds to $f$. 

From Proposition \ref{extension_of_b-surface}, we also have
\begin{Prop} \label{Prop:conjecture_holds}
  Conjecture \ref{Conj:b-surface} holds for the singular self-dual Zollfrei conformal 
  structure corresponding to $f\in\cS(\R^2)^{\text{\rm sym}}$. 
\end{Prop}

\section{Construction of singular twistor spaces} \label{Section:twistor}

Recall that we identify $\RP^3-\{y_0\}$ with the line bundle 
$\cO_\R(1)\rightarrow\RP^2$, where $y_0=[0:0:0:1]$. 
Let $\cS(\R)^{\text{odd}}$ be the subset of $\cS(\R)$ consisting of 
odd functions. 
For each $s(t)\in\cS(\R)^{\text{odd}}$, we define a smooth section 
$\tilde{s}$ of the line bundle $\cO_\R(1)\rightarrow\RP^2$ by 
$\tilde{s}\([0:0:1]\) = [0:0:1:0] $ and
\begin{equation} \label{eq:section}
 \tilde{s}\(\[-\sin\chalf : \cos\chalf :t \]\)
   = \[-\sin\chalf : \cos\chalf : t : s(t)\], 
\end{equation}
Then we put 
 $$ P=\cO_\R(1)+ i\, \tilde{s}(\RP^2)
   = \left\{u+i\, \tilde{s}(x)\in \cO(1) : x\in \RP^2, u\in\cO_\R(1)_x  \right\}, $$
which is a deformation of $\cO_\R(1)$ in $\cO(1)$. 
Adding the point $y_0$ to the pair $(\cO(1),P)$, we have the 
twistor space $(\CP^3,\hat{P})$. 
Notice that $\hat{P}$ is the image of a map $\iota : \RP^3\rightarrow\CP^3$ 
given by $\iota(y_0)=y_0$ and 
$$ \cO_\R(1) \ni u \longmapsto u+i\, \tilde{s}(\pi(u)) \in P $$
using the projection $\pi$. 
For the later convenience, we use the pure imaginary valued function 
$h(t)=i s(t)\in i \cS(\R)^{\text{odd}}$ from now on. 
Then our goal is : 
\begin{Thm}\label{Thm:twistor_correspondence}
  There is a natural one-to-one correspondence between 
  \begin{itemize}
   \item singular self-dual Zollfrei conformal structures on $\widetilde{Gr}_2(\R^4)$ 
      corresponding to $f(x)\in \cS(\R^2)^{\text{\rm sym}}$, 
	  \item the pairs $(\CP^3,\hat{P})$ corresponding to 
   $h(t)\in i\cS(\R)^{\text{\rm odd}}$,  
  \end{itemize}
  which satisfies the following property. 
  There is a double fibration $\widetilde{Gr}_2(\R^4)\overset{\wp}
  \longleftarrow \hat{\cZ} \overset{\Psi}\longrightarrow \CP^3$ 
  such that 
 \begin{enumerate}
  \item $\Psi$ is a continuous surjection and $\wp$ is a complex disk bundle, 
   i.e. for each $x\in\widetilde{Gr}_2(\R^4)$, $\hat{\cZ}_x=\wp^{-1}(x)$ is 
   biholomorphic to the complex unit disk, 
  \item $\Psi_x : \hat{\cZ}_x \rightarrow\CP^3$ is holomorphic on the interior of 
   $\hat{\cZ}_x$ and $\Psi_x(\del \hat{\cZ}_x)\subset \hat{P}$, 
   where $\Psi_x$ is the restriction of $\Psi$ on $\hat{\cZ}_x$,
  \item the restriction of $\Psi$ on $\hat{\cZ}-\del\hat{\cZ}$ is bijective 
        onto $\CP^3-\hat{P}$, 
  \item $\{\wp(\Psi^{-1}(y))\}_{y\in P}$ is equal to the set of $\b$-surfaces 
        on $\widetilde{Gr}_2(\R^4)$. 
 \end{enumerate}  
  Moreover, this correspondence is explicitly given by the formulas 
  $f=2i \(\pd{t}h\)^\vee$ and 
  $h=-\frac{1}{4}\cH \hat{f}$.  
\end{Thm}

\noindent
{\it Remark :} The dual Radon transform $(\cdot)^\vee$ and the 
Hilbert transform $\cH$ are explained in Section \ref{Appendix}. 

\vspace{1ex}
Here we explain some reasons why the above construction of the 
twistor space is reasonable. 
First recall that our singular metrics are standard on the equator $W_0$, 
so the twistor spaces would be standard over the corresponding point 
$[0:0:1] \in\RP^2$. This corresponds to $\tilde{s}([0:0:1])=[0:0:1:0]$. 
Moreover, our singular metrics have $S^1$-symmetry, so the twistor spaces 
would also have a similar symmetry. 
From the twistor correspondence for the standard Zoll projective structure, 
the $S^1$-action on $\RP^2$ is given by
$$ S^1\ni\t : [z_1:z_2:z_3] \mapsto 
[\cos\t\ z_1-\sin\t\ z_2 : \sin\t\ z_1+\cos\t\ z_2 : z_3]. $$
And the lift of this $S^1$-action on 
$\cO_\R(1)=\RP^3-\{y_0\}$ is given by
$$  [z_1:z_2:z_3:z_4] \mapsto 
[\cos\t\ z_1-\sin\t\ z_2 : \sin\t\ z_1+\cos\t\ z_2 : z_3 : z_4]. $$
So $\tilde{s}$ should be $S^1$-equivariant, if the twistor space is given by the section $\tilde{s}$ of $\cO_\R(1)$, and if it corresponds to our singular metric. 
Such a section is given by an odd function $s(t)$ 
in the manner of (\ref{eq:section}). \\

Before we start to prove Theorem \ref{Thm:twistor_correspondence}, 
we make some remarks here. 
First, we study about the real twistor correspondence for our 
singular self-dual Zollfrei metric. 
This real correspondence has important meaning not only for its geometric
significance but also as a step of the construction of the complex correspondence.
We define a map $\Phi : F\rightarrow \cO_\R(1)$ by 
\begin{equation}\label{deformed_Phi}
 \begin{aligned}
     F|_{D_\pm}\ni \(x_1^\pm,x_2^\pm,x_3^\pm,x_4^\pm,\z\) & \longmapsto 
      \[ 1 : \z : -x_1^\pm -x_2^\pm \z : 
         -x_3^\pm \z + x_4^\pm \pm\frac{\vp}{\sin\s} \], \\
    F|_{W_0} \ni \(\a,0,\e_1,\e_2,\x\) & \longmapsto \[ \rule{0in}{2.5eX}
      -\sin\a : \cos\a : \x : -\e_1\x+\e_2+\frac{1}{4}\hat{f}(\x) \], 
 \end{aligned}
\end{equation}
where $\z=-\cot\s$ and $W_0=\{\b=0\}\subset W$. 
From (\ref{eq:beta-surfaces}) and (\ref{coordinate_changed_eq}), 
the inverse image of a point by this map is a $\b$-surface. 
Hence $\cO_\R(1)$ is identified with the space of 
$\b$-surfaces by this map. 

Next, we study about the complex structure on $\cZ$, 
where $\cZ$ is given in the diagram (\ref{the_diagram}). 
$\cZ$ has a 
natural complex structure defined from the self-dual metric on $TS^2$ 
as the following (cf.\cite{bib:LM05}, see also \cite{bib:BDM}). 
Recall that for a given Petean's metric, any $\b$-plane at a point 
is written in the form $\langle \Gn_1,\Gn_2 \rangle$ where $\Gn_1,\Gn_2$ 
is given by (\ref{n_1n_2}). 
Putting $\z=-\cot \s$ and by the ``analytic continuation", 
we obtain the complex tangent vectors
$$ \Gn'_1 = -\z^+\pd{x_1^+}+\pd{x_2^+}
       -\frac{f}{2} \( -\z^+\pd{x_3^+}+\pd{x_4^+} \), \quad 
  \Gn'_2 = \pd{x_3^+}+\z^+\pd{x_4^+}, $$
on $\cZ |_{TD_+}$. The complex structure on $\cZ |_{TD_+}$ is 
defined so that its (0,1)-vectors are given by 
$\Gn'_1, \Gn'_2$ and $\pd{\bar{\z^+}}$. 

The key to prove Theorem \ref{Thm:twistor_correspondence} 
is to construct a map $\Phi_\C: \cZ\rightarrow \cO(1)$ 
which is described in the form 
$$ \(x_1^+,x_2^+,x_3^+,x_4^+,\z^+\) \longmapsto 
      \[ 1 : \z^+ : -x_1^+ -x_2^+ \z^+ : 
         -x_3^+ \z^+ + x_4^+ +H(x_1^+,x_2^+,\z^+) \] $$ 
on $D_+$. 
We will construct the surjection $\Psi:\hat{\cZ}\rightarrow\CP^3$ 
as an extension of $\Phi_\C$. 
Notice that $\Phi_\C$ is holomorphic on $\cZ |_{D_+}$ with respect to 
the above complex structure if and only if  
(i) $H(x_1^+,x_2^+,\z^+)$ is holomorphic for $\z^+$ and (ii) $H$ solves
\begin{equation} \label{eq_of_holomorphicity}
 \(-\z^+\pd{x_1^+}+\pd{x_2^+}\) H(x_1^+,x_2^+,\z^+)
    = \frac{f(x_1^+,x_2^+)}{2} \, \((\z^+)^2+1\). 
\end{equation}

\section{Holomorphic disks} \label{Section:holo_disks}

We prove the following in this section. 
\begin{Prop} \label{Prop:holo_disks}
 The pair $(\CP^3,\hat{P})$ corresponding to a function 
 $h(t)\in i\cS(\R)^{\text{\rm odd}}$ satisfies the condition $(\sharp)$. 
\end{Prop}

We first study about the holomorphic disks which 
do not contain the singular point $y_0$. 
Recall that any holomorphic disk in $(\CP^2,\RP^2)$ whose relative homology 
class generates $H_2(\CP^2,\RP^2;\Z)\cong\Z$ is given by 
one of the hemispheres of a complex line (p.498-500 of \cite{bib:LM02}). 
We call such a disk the {\it standard disk} in $(\CP^2,\RP^2)$. 
Remark that $\hat{P}$ is homotopic to the standard $\RP^3$ in $\CP^3$, 
so $H_2(\CP^3,\hat{P}; \Z) \simeq \Z$. 

\begin{Lem} \label{Lem:holo_lift}
 Let $D$ be the complex unit disk and 
 $\vp: D \rightarrow \cO(1)$ be a continuous map with $\vp(\del D)\subset P$. 
 If $\vp$ is holomorphic on the interior of $D$ and the relative homology class 
 $[\vp]$ generates $H_2(\CP^3,\hat{P};\Z)$, then 
 $\vp$ is a holomorphic lift of some standard disk in $(\CP^2,\RP^2)$, i.e. 
 the image of the composition $\pi\action\vp$ is a standard disk where 
 $\pi:\cO(1)\rightarrow\CP^2$ is the projection. 
\end{Lem}
\begin{proof}
 Let $i:\cO(1)\rightarrow\CP^3$ be the inclusion, then we have the isomorphisms 
 $$ \Z\cong H_2(\CP^2,\RP^2) \overset{\pi_*}\longleftarrow 
      H_2(\cO(1),P) \overset{i_*}\longrightarrow H_2(\CP^3,\hat{P})\cong\Z.$$
 Since $[i\action\vp]$ generates $H_2(\CP^3,\hat{P})$, $[\pi\action\vp]$ 
 also generates $H_2(\CP^2,\RP^2)$. 
 Because $\pi\action\vp$ defines a holomorphic disk in $(\CP^2,\RP^2)$, 
 it is a standard disk. 
\end{proof}

Now let $D$ be a standard disk in $(\CP^2,\RP^2)$ and $\cF_D$ be the set of 
holomorphic lifts of $D$ in $(\cO(1),P)$. In the proof of next Lemma, 
we follow to the method of LeBrun and Mason
(\cite{bib:LM02},\cite{bib:LM05}). 

\begin{Lem} \label{Lem:R^2_families}
 $\cF_D$ has a structure of smooth family of holomorphic disks 
 parameterized by $\R^2$. 
\end{Lem}

\noindent
{\it Remark :}
In our situation, we can calculate explicitly, and we have no need to use the 
inverse function theorem of Banach space. 
So we can treat many metrics which are far from the standard one. 

\begin{proof}[Proof of \ref{Lem:R^2_families}]
Recall that $\cO(1)\cong\CP^3-\{y_0\}$ \ ($y_0=[0:0:0:1]$) and 
the projection $\pi:\cO(1)\rightarrow\CP^2$ is given by 
$[z_1:z_2:z_3:z_4]\mapsto[z_1:z_2:z_3]$. 
Let $U$ be the affine open set in $\cO(1)$ 
whose coordinate is defined by
$$ \Gz_1=\frac{z_2-iz_1}{z_2+iz_1}, \quad 
   \Gz_2=\frac{z_3}{z_2+iz_1}, \quad
   \Gz_3=\frac{z_4}{z_2+iz_1}. $$
Then the intersection $B:=\cO_\R(1)\cap U$ 
is equal to the set given by 
$$ \Gz_1\bar{\Gz_1}=1,\quad \Gz_1\bar{\Gz_2}=\Gz_2, \quad 
   \Gz_1\bar{\Gz_3}=\Gz_3. $$
We can parameterize $B$ by $(\Gz_1,\Gz_2,\Gz_3)=
(e^{i\c},t_1e^{\frac{i\c}{2}},t_2e^{\frac{i\c}{2}})$ by using 
$\c$, $t_1$ and $t_2$, 
or equivalently $z=[-\sin\chalf : \cos\chalf : t_1 : t_2]$. 
So we have $B\cong \R^3/\Z$ where the $\Z$-action is generated by
$(\c,t_1,t_2)\mapsto(\c+2\pi,-t_1,-t_2)$. 
In the similar way, we can parameterize $B':=P\cap U$ by 
$$ (\Gz_1,\Gz_2,\Gz_3)=
(e^{i\c},t_1e^{\frac{i\c}{2}},[t_2+h(t_1)]e^{\frac{i\c}{2}}). $$ 
Notice that $(\Gz_1,\Gz_2)$ defines a coordinate on the open set
$\pi(U)$ in $\CP^2$, and that 
$\pi(B')=\pi(B)\cong \R^2/\Z$ is equal to $\RP^2-\{[0:0:1]\}$. 

The boundary $\del D$ of a standard holomorphic disk $D$ is a real projective line 
in $\RP^2$. 
Each real line in $\RP^2$ is described by either form of 
\begin{enumerate}
 \item $(\Gz_1,\Gz_2)=(e^{i\c},a+\bar{a}e^{i\c})$ \ for some $a\in\C$, 
 \item $(\Gz_1,\Gz_2)=(e^{2i\a},\x e^{i\a})$ \ for some $\a\in S^1=\R/2\pi\Z$, 
\end{enumerate}
where $\c\in S^1$ and $\x\in\R$ are the parameters. 
The case (2) occurs when the circle passes through $[0:0:1]$. 
Notice that in the case (2), $\a$ and $\a+\pi$ correspond to the same 
line with opposite orientations. \\

We give the proof for the two cases separately; 
one is the case when $\del D$ is given by (1) and the other case is by (2).  
In the case (1), any lift of $\del D$ is described by 
\begin{equation}\label{Gz's_case_1} 
  (\Gz_1,\Gz_2,\Gz_3)=\(e^{i\c},a+\bar{a}e^{i\c},
  \[u(\c)+h\(ae^{-\frac{i\c}{2}}+\bar{a}e^{\frac{i\c}{2}}\)\]
  e^{\frac{i\c}{2}} \) 
\end{equation}
with the parameter $\c\in S^1$, 
where $u(\c)$ is an unknown real function satisfying $u(\c+2\pi)=-u(\c)$. 
Notice that for a fixed $a\in\C$, there  two possibilities of $D$, 
i.e. the upper and lower hemispheres. 
Corresponding to this, we want to choose $u(\c)$ so that 
the circle (\ref{Gz's_case_1}) extends 
holomorphically to the interior or exterior region of $\{|\w|=1\}$ where 
$\w=e^{i\c}$.
In the interior case,  we want to choose $u(\c)$ such that
$\Gz_3$ contains no negative power in its Fourier expansion. 
Because $u$ is real valued and $h$ is pure imaginary valued, we can expand  
\begin{equation} \hspace{15mm}
  u(\c) =\sum_{l=0}^\infty \left\{ u_l\ e^{i\c(l+\frac{1}{2})}+
            \bar{u}_l\ e^{-i\c(l+\frac{1}{2})}\right\}, 
\end{equation} \vspace{-2mm}
\begin{equation} \label{expand_h}
  h\(ae^{-\frac{i\c}{2}}+\bar{a}e^{\frac{i\c}{2}}\)
           = \sum_{l=0}^\infty \left\{ h_{a,l}\ e^{i\c(l+\frac{1}{2})}-
            \bar{h}_{a,l}\ e^{-i\c(l+\frac{1}{2})} \right\}.
\end{equation}
So if we put $u_l=h_{a,l}$ for $l \geq 1$, then we have 
$$ \Gz_3=
  \[u(\c)+h\(ae^{-\frac{i\c}{2}}+\bar{a}e^{\frac{i\c}{2}}\)\] e^{\frac{i\c}{2}} 
  =2\sum_{l=0}^{\infty} h_{a,l}\ e^{i\c(l+1)} + \k + \bar{\k} e^{i\c}, $$
where $\k=\bar{u}_0-\bar{h}_{a,l}$ can be taken as an arbitrary 
complex constant. 
The circle (\ref{Gz's_case_1}) is described in the homogeneous coordinate by 
\begin{equation} \label{holo_disk_interior}
  [1:\z: -2(\Imag a)+2(\Real a)\z : -2(\Imag \k)+2(\Real \k)\z +F(a,\w(\z))] 
\end{equation}
with the parameter $\z\in\R\cup\{\infty\}$, where 
\begin{equation} \label{z&w}
 \w(\z)=\frac{\z-i}{\z+i}, \qquad
 F(a,\w)=\frac{4i}{1-\w}\sum_{l=0}^\infty h_{a,l}\ \w^{l+1}.
\end{equation}
The circle (\ref{holo_disk_interior}) extends holomorphically to 
$\{\Imag \z\geq 0\}$ and $\{|\w|\leq 1\}$, 
and hence $\cF_D$ is a smooth family parameterized by $\k\in\C\simeq\R^2$.  
 
The exterior case is, in the same way, 
given by (\ref{holo_disk_interior}) using $\overline{F(a,\bar{\w}^{-1})}$ 
instead of $F(a,\w)$. 
In this case, (\ref{holo_disk_interior}) 
extends to $\{\Imag \z\leq 0\}$ and $\{|\w|\geq 1\}$. \\

In the case (2), it is enough to consider the standard holomorphic 
disk given by $\{\x\in\C : \Imag \x\geq 0\}$. 
In this case, any lift of $\del D$ is described by 
\begin{equation} \label{holo_disk_W0}
 (\Gz_1,\Gz_2,\Gz_3)=(e^{2i\a}, \x e^{i\a}, [u(\x)+h(\x)]e^{i\a}) 
\end{equation}
with the parameter $\x\in\R$, where $u(\x)$ is an unknown real function.  
We define a holomorphic function $G(\x)$ on
$\{\Imag\x >0\}$ by
$$ G(\x)=\frac{1}{\pi i} \int_{-\infty}^{\infty} \frac{h(\mu)}{\mu-\x}d\mu. $$
Then $G(\x)$ extends smoothly to $\{\Imag\x =0\}$ by 
\begin{equation} \label{G_on_the_real}
  G(\x)=\frac{1}{\pi i}\ \pv\int_{-\infty}^{\infty} \frac{h(\mu)}{\mu-\x} d\mu
  + h(\x). 
\end{equation}
The first term of the right hand side is equal to $-\cH h(\x)$ 
where $\cH$ is the Hilbert transform defined by (\ref{Hilbert_transf}). 
If we put $v(\x)=u(\x)+\cH h(\x)$, then we have 
$$ u(\x)+h(\x)= v(\x)+G(\x) $$ 
on $\{\Imag\x =0\}$. The circle (\ref{holo_disk_W0}) is described 
in the homogeneous coordinate by 
\begin{equation}\label{holo_disk_case2}
 [-\sin\a: \cos\a : \x : v(\x)+G(\x)] 
\end{equation}
with the parameter $\x\in\R$. 
If the circle (\ref{holo_disk_case2}) is bounded by a holomorphic disk, 
$v(\x)$ should extend holomorphically to the upper half plane of $\x$. 
Because $v(\x)$ is real valued on $\R$, it extends to a 
holomorphic function on whole $\C$. Hence $v(\x)$ expands to 
a non-negative power series of $\x$ with real coefficients. 

Moreover (\ref{holo_disk_case2}) converges to 
$y_0=[0:0:0:1]$ by $\x\rightarrow \infty$, 
if and only if $v(\x)$ has any term higher than or equal to $\x^2$. 
These cases are removable, because the disk is not contained in $\cO(1)$.   
Hence $v(\x)$ should be a degree-one polynomial. 
Because the ambiguity is given by the coefficients of $v(\x)$, 
$\cF_D$ is a smooth $\R^2$-family of holomorphic disks. 
\end{proof}

Let $\cF$ be the set of all the holomorphic disks in $(\cO(1),P)$ 
such that the relative homology class of each disk 
generates $H_2(\CP^3,\hat{P};\Z)$. 
Then we have $\cF=\cup \cF_D$ from Lemma \ref{Lem:holo_lift}. 

\begin{Lem} \label{Lem:para_by_TS^2}
 $\cF$ has a structure of smooth family of holomorphic disks 
 parametrised by $TS^2$ such that interiors of the disks of $\cF$ foliate 
 $\cO(1)-\cO(1)|_{\RP^2}$. 
\end{Lem}
\begin{proof}
Recall the diagram (\ref{the_diagram}), and we construct 
a smooth map $\Phi_\C : \cZ\rightarrow \cO(1)$ as 
a deformation of $\Phi_{\C,0}$ such that 
$\{\Phi_\C(\tilde{\Gp}^{-1}(x))\}_{x\in TS^2}$ is equal to $\cF$. 
Using the coordinate defined in Section \ref{Section:Standard},
we define $\Phi_\C$ on $\cZ|_{D_+}$ by 
\vspace{2mm}
\begin{equation} \label{Phi_TD+}
 \(x_1^+,x_2^+,x_3^+,x_4^+,\z^+\) \longmapsto 
   \[ 1 : \z^+ : -x_1^+ -x_2^+ \z^+ : 
        -x_3^+ \z^+ +x_4^+ +H^+(x_1^+,x_2^+,\z^+)\], 
\end{equation}
$$ H^+(x_1^+,x_2^+,\z^+) = F\(a,\w(\z^+)\), 
  \quad a=\frac{i}{2}(x_1^+ +i x_2^+),$$
where $F$ is given in (\ref{z&w}). 
Similarly we define $\Phi_\C$ on $\cZ|_{D_-}$ by the same formula as 
(\ref{Phi_TD+}) by exchanging the sign $\cdot^+$ by $\cdot^-$ and 
$$ H^-(x_1^-,x_2^-,\z^-) = -\overline{F\(a,\bar{\w(\z^-)}^{-1}\)}, \quad
   a=\frac{i}{2}(x_1^- +i x_2^-). $$
Remark that the parameters run $\{\Imag\z^+\geq 0\}$ and 
$\{\Imag\z^-\leq 0\}$. 
On the other hand, 
we define $\Phi_\C$ on $\cZ|_{W_0}=\{\b=0\}$ by 
\begin{equation} \label{Phi_W0}
 \(\a,0,\e_1,\e_2,\x\) \longmapsto 
 \[ \hasira -\sin\a : \cos\a : \x : -\e_1\x+\e_2 + G(\x) \]. 
\end{equation}
From (\ref{holo_disk_interior}) and (\ref{holo_disk_case2}), we obtain that 
$\{\Phi_\C(\tilde{\Gp}^{-1}(x))\}_{x\in TS^2}$ is equal to $\cF$. 
Hence $\cF$ is parametrised by $TS^2$. 

We have to prove that the above $\Phi_\C$ is smooth. 
We now check that (\ref{Phi_TD+}) and (\ref{Phi_W0}) are continued. 
Here we omit the sign `+' for the simplicity. 
Using the coordinate change 
$$ \x=-i\cot\b\ \frac{\w-e^{2i\a}}{\w+e^{2i\a}}, $$
which is obtained from  (\ref{coordinate_change}) and (\ref{z&w}), and so on, 
we have
$$ \begin{aligned}
   \[ 1 : \z \right.&\left. : -x_1-x_2\z :  -x_3\z+x_4+H(x_1,x_2,\z)\] \\
  =& \[ \x\cos\a\tan\b-\sin\a : -\x\sin\a\tan\b +\cos\a : \x :
        -\e_1\x+\e_2 + B(\a,\b,\x) \hasira \],
  \end{aligned} 
$$ 
where
\begin{equation}\label{eq:B}
 B(\a,\b,\x) = \frac{4e^{i\a}}{\w+e^{2i\a}}\sum_{l=0}^\infty 
 h_{a,l}\ \w^{l+1}. 
\end{equation}
In general, we define the operator $\Pi$ on the $L^2$-functions 
on $S^1=\{|\w|=1\}$ by 
$$\Pi : \ u(\w)= \sum_{k=-\infty}^\infty u_k \w^k
  \ \longmapsto\ \sum_{k=0}^\infty u_k \w^k.$$
As explained in \cite{bib:Taylor}, when $|\w|<1$ we have
$$ \Pi u(\w)=\frac{1}{2\pi i} \int_{S^1} \frac{u(\y)}{\y-\w}\ d\y 
  = \frac{1}{2\pi} \int_{S^1} \frac{u(e^{i\c})}{e^{i\c}-\w}\ e^{i\c}d\c. $$
Hence, in our case, 
$$ \begin{aligned}
  \sum_{l=0}^\infty h_{a,l}\ \w^l 
  &= \Pi\(e^{-\frac{i\c}{2}}h\(ae^{-\frac{i\c}{2}}+\bar{a}e^{\frac{i\c}{2}}\)\) 
       (\w) \\
  &= \frac{1}{2\pi}\int_{S^1} 
  \frac{e^{\frac{i\c}{2}}h\(ae^{-\frac{i\c}{2}}+\bar{a}e^{\frac{i\c}{2}}\)}
  {e^{i\c}-\w} \ d\c. 
 \end{aligned} $$
If we change the parameter $\c$ by
$\mu=ae^{-\frac{i\c}{2}}+\bar{a}e^{\frac{i\c}{2}}
 =\cot\b\sin(\chalf-\a)$, then we have
 $$ \begin{aligned}
   B(\a,\b,\x)
    &= \frac{2e^{i\a}\w}{\pi(\w+e^{2i\a})} \int_{-\cot\b}^{\cot\b}
	      \frac{h(\mu)}{e^{\frac{i\c}{2}}-\w e^{-\frac{i\c}{2}}} \cdot
		  \frac{2d\mu}{\cot\b\ \cos(\chalf-\a)} \\
    &= \frac{1+\x^2\tan^2\b}{\pi i} \int_{-\cot\b}^{\cot\b}
	      \frac{h(\mu)}{\mu-\x\sqrt{1-\mu^2\tan^2\b}} \cdot
		  \frac{d\mu}{\sqrt{1-\mu^2\tan^2\b}}. 
  \end{aligned} $$
Because $h(t)$ is rapidly decreasing, $B(\a,\b,\x)$ extends continuously to 
$\b=0$ and we have $B(\a,0,\x)=G(\x)$. 
Hence $\Phi_\C$ is continuous on $\cZ|_{D_+\cup W_0}$. 
In the similar argument, $\Phi_\C$ is continuous on $\cZ|_{D_-\cup W_0}$. 
Moreover it is smooth because of the above formula. 

We check that $\cO(1)-\cO(1)|_{\RP^2}$ is foliated by $\cF$. 
For any $u\in\cO(1)-\cO(1)|_{\RP^2}$, there is a unique standard 
disk $D$ in $(\CP^2,\RP^2)$ which contains $\pi(u)\in\CP^2-\RP^2$. 
Then, from (\ref{holo_disk_interior}) or (\ref{holo_disk_case2}), there 
is a unique holomorphic disk in $\cF_D$ which contains $u$. 
\end{proof}

\noindent
{\it Remark : }
The smoothness of $\Phi_\C$ is also checked in the proof of Lemma \ref{Lem:h=>f}, 
actually $\Phi_\C$ is holomorphic on the interior of $\cZ$ 
for some complex structure on there. 

\noindent
{\it Remark : }
 For each $t\in S^2$, the set of holomorphic disks 
 $\{\Phi_\C(\tilde{\Gp}^{-1}(x))\}_{x\in T_tS^2}$ 
 is equal to $\cF_D$, where $D=\Gq(\Gp^{-1}(t))$ is the standard disk 
 corresponds to $t$ in the twistor correspondence for the standard Zoll projective 
 structure. 

\begin{proof}[Proof of \ref{Prop:holo_disks}]
Notice that 
$\CP^3-\hat{P}=(\cO(1)-\cO(1)|_{\RP^2}) \coprod (\cO(1)|_{\RP^2}-P)$. 
We require a family of holomorphic disks through $y_0$ whose interiors 
foliate $\cO(1)|_{\RP^2}-P$. 
Similar to the standard case, such a family is given by the disks 
\begin{equation} \label{infinitely_far_disks}
 \{ [z_1:z_2:z_3:u+h(z_3)] : \Imag u\geq 0]\} \cup \{y_0\} 
\end{equation}
for some $(z_1,z_2,z_3)\in\R^3-\{0\}$. 
There is a one-to-one correspondence between $S^2_\infty$ and the above 
holomorphic disks, which is given by putting $(z_1,z_2,z_3)=t\times v$
for each point in $S^2_\infty$ of the form (\ref{pt_in_S^2_infty}). 

If we define $\hat{\cF}$ as the union of $\cF$ and the disks of the form 
(\ref{infinitely_far_disks}), 
then $\hat{\cF}$ is a continuous family of holomorphic disks in $(\CP^3,\hat{P})$ 
parameterized by $\widetilde{Gr}_2(\R^4)$. 
The conditions ($\sharp\,1$), ($\sharp\,2$), and the uniqueness for $\hat{\cF}$ 
follow directly from the construction. 
\end{proof}

\section{Twistor correspondence} \label{Section:corresp}
In this section, we give the proof of Theorem \ref{Thm:twistor_correspondence}. 

\begin{Lem} \label{Lem:h=>f}
 For a given twistor space $(\cO(1),P)$ corresponding to a function 
 $h\in i\cS(\R)^{\text{\rm odd}}$, 
 there are a smooth map $\Phi_\C: \cZ\rightarrow\cO(1)$ and 
 a function $f\in\cS(\R^2)^{\text{\rm sym}}$ which satisfy 
 the following properties: 
 \begin{enumerate}
  \item $\Phi_{\C,x} : \cZ_x\rightarrow\cO(1)$ is holomorphic on the interior 
        of $\cZ_x=\tilde{\Gp}^{-1}(x)$ and $\Phi_{\C,x} (\del\cZ_x)\subset P$, 
		where $\Phi_{\C,x}$ is the restriction of $\Phi_\C$ on $\cZ_x$, 
  \item $\Phi_\C$ is injective on $\cZ-\del\cZ$, 
  \item $\{\tilde{\Gp}(\Phi_\C^{-1}(y))\}_{y\in P}$ is equal to 
        the set of $\b$-surfaces on $TS^2$, 
 \end{enumerate} 
 respecting the self-dual metric on $TS^2$ and the complex structure on $\cZ$ 
 corresponding to $f$. 
 Such $f$ is given by 
 \begin{equation} f(x)=2i\(\PD{h}{t}\)^\vee (x). \end{equation}
\end{Lem}
\begin{proof}
The map $\Phi_\C$ is already constructed in Lemma \ref{Lem:para_by_TS^2}, 
and the conditions (1) and (2) are already checked. 
We now construct $f\in\cS(\R^2)^{\text{sym}}$ and check the 
condition(3). 
If $\Phi_\C$ is holomorphic on the interior of $\cZ$ with respect to the 
complex structure defined from $f$, then the equation 
(\ref{eq_of_holomorphicity}) holds, and this is equivalent to 
\begin{equation} \label{holomorphicity_F}
 \(\w\pd{a}-\pd{\bar{a}}\) \(\hasira (1-\w)F(a,\w)\) = -2\w f. 
\end{equation}
For given $h$, the function $f$ is uniquely defined by (\ref{holomorphicity_F}). 
Actually, from the identity 
$$ \(e^{\frac{i\c}{2}}\pd{a}-e^{-\frac{i\c}{2}}\pd{\bar{a}}\)
  h\(ae^{-\frac{i\c}{2}}+\bar{a}e^{\frac{i\c}{2}}\)=0, $$ 
we obtain
$$ \PD{h_{a,l}}{\bar{a}} = \PD{h_{a,l-1}}{a}  \quad l\geq 1,
   \ \text{and} \quad
   \PD{h_{a,0}}{\bar{a}} \ \text{is real valued}. $$
So the equation (\ref{holomorphicity_F}) holds if and only if 
$f(x)=2i\PD{h_{a,0}}{\bar{a}}$, 
where $x=(x_1^+,x_2^+)$ and $a=\frac{i}{2}(x_1^++ix_2^+)$. 
Here we have 
$$ \begin{aligned}
  \PD{h_{a,0}}{\bar{a}}
   &= \frac{1}{2\pi} \pd{\bar{a}}\int_{S^1} e^{-\frac{i\c}{2}} 
		  h\(ae^{\frac{-i\c}{2}}+\bar{a}e^{\frac{i\c}{2}}\) \ d\c \\[0.5ex]
   &= \frac{1}{2\pi} \int_{S^1} \PD{h}{t}
            \(ae^{\frac{-i\c}{2}}+\bar{a}e^{\frac{i\c}{2}}\) \ d\c \\
   &= \(\PD{h}{t}\)^\vee (x), 
 \end{aligned} $$
hence we put 
\begin{equation}
  f(x)=2i\(\PD{h}{t}\)^\vee (x) 
\end{equation}
which is real valued, rapidly decreasing and axisymmetric. 

Now we prove that the condition (3) holds for this $f$. 
The equation (\ref{eq_of_holomorphicity}) holds on 
$\{\z^+\in\C : \Imag \z^+=0\}$, so, for a fixed $\z^+\in\R$, the functions 
$$ -x_1^+-x_2^+\z^+ \quad \text{and} \quad 
   -x_3^+\z^+ +x_4^+ + H(x_1^+,x_2^+,\z^+) $$
are annihilated by $\Gn_1,\Gn_2$ given by (\ref{n_1n_2}). 
Hence $\tilde{\Gp}(\Phi^{-1}(y))$ is a $\b$-surface on $TS^2|_{D_+}$. 
In the same way, $\tilde{\Gp}(\Phi^{-1}(y))$ is a $\b$-surface on $TS^2|_{D_-}$ 
so $\tilde{\Gp}(\Phi^{-1}(y))$ is a $\b$-surface on $TS^2$. 
Hence the condition (3) holds. 
\end{proof}

\begin{Lem} \label{Lem:f=>h}
 For a given self-dual metric on $TS^2$ corresponding to a function 
 $f\in\cS(\R^2)^{\text{\rm sym}}$, 
 there are a function $h\in i\cS(\R)^{\text{\rm odd}}$ and 
 a smooth surjection $\Phi_\C :\cZ\rightarrow\cO(1)$ 
 which satisfy the condition {\rm (1), (2)} and {\rm (3)} in the 
 Lemma \ref{Lem:h=>f} respecting the complex structure 
 on $\cZ$ defined from the self-dual metric and the twistor space $(\cO(1),P)$ 
 corresponding to $h$. 
 Such $h$ is given by 
  \begin{equation} \label{h_from_f}
  h(t)=-\frac{1}{4} \cH \hat{f}(t). 
 \end{equation}
\end{Lem}
\begin{proof}
We have already seen in (\ref{deformed_Phi}) that there is a natural map 
$\Phi:F=\del\cZ\rightarrow\cO_\R(1)$, which is given by 
\begin{equation} \label{Phi_on_W_0}
 (\a,0,\e_1,\e_2,\x) \longmapsto 
 \[-\sin\a : \cos\a : \x :-\e_1\x+\e_2+\frac{1}{4}\hat{f}(\x) \],
\end{equation}
on $F|_{W_0}$. 
Now we deform this map so as to extend holomorphically to 
the upper half plane of $\x$. For this purpose, we put
\begin{equation} \label{f=>h}
  h(\x)=\frac{1}{4\pi i} \pv \int_{-\infty}^\infty
    \frac{\hat{f}(\mu)}{\mu-\x} \ d\mu
	    = -\frac{1}{4} \cH \hat{f}(\x) 
\end{equation}
which is odd, rapidly decreasing and pure imaginary valued. 
Similar to the previous section, 
$G(\x)=h(\x)+\frac{1}{4}\hat{f}(\x)$ 
extends holomorphically to $\{\Imag \x> 0\}$ by
$$ G(\x)=\frac{1}{4\pi i} \int_{-\infty}^\infty 
\frac{\hat{f}(\mu)}{\mu-\x} d\mu. $$
Instead of (\ref{Phi_on_W_0}), we define a map 
$\Upsilon : \cZ|_{W_0}\rightarrow\cO(1)$ by 
$$ (\a,0,\e_1,\e_2,\x) \longmapsto 
  \[ -\sin\a : \cos\a : \x : -\e_1\x+\e_2+G(\x) \]. $$
Let $(\cO(1),P)$ be the twistor space corresponding to 
$h\in i\cS(\R)^{\text{odd}}$, 
then $\Upsilon_x: \cZ_x\rightarrow\cO(1)$ is a  
holomorphic disk in $\cO(1)$ with boundary on $P$ for each $x\in TS^2|_{W_0}$. 

By Lemma \ref{Lem:h=>f}, we can construct a smooth map 
$\Phi_\C : \cZ \rightarrow \cO(1)$ with $\Phi_\C(\partial\cZ)=P$, 
which is holomorphic on the interior of $\cZ$ with respect to the 
complex structure defined from the function
$$ k(x)=2i\(\PD{h}{t}\)^\vee (x). $$
Connecting with (\ref{f=>h}), we have
$$ k(x)=\frac{1}{2i}\(\pd{t}\cH \hat{f} \)^\vee (x). $$
By the inversion formula (\ref{inversion_formula}), we have $k(x)=f(x)$. 
Therefore $\Phi_\C$ is holomorphic with respect to the 
given complex structure on $\cZ$. 
Remark that the restriction of $\Phi_\C$ to $\cZ|_{W_0}$ 
is equal to $\Upsilon$. 
The conditions (1), (2) and (3) follow from Lemma \ref{Lem:h=>f}. 
\end{proof}

\begin{Lem} \label{Lem:one-to-one}
 (\ref{h_from_f}) defines a one-to-one correspondence between 
 $f(x)\in\cS(\R^2)^{\text{\rm sym}}$ 
 and $h(\x)\in i\cS(\R)^{\text{\rm odd}}$. 
\end{Lem}
\begin{proof}
 From the Schwartz's theorem (Theorem 2.4 of \cite{bib:Helgason}), 
 we can check that the Radon transform 
 $f(x) \rightarrow \vp(p)=\hat{f}(p)$ is one-to-one from 
 $\cS(\R^2)^{\text{\rm sym}}$ to $\cS(\cP)^{\rm sym}$, 
 where $\cP$ is the set of unoriented lines in $\R^2$ and 
 $\cS(\cP)^{\rm sym}$ is the set of rapidly decreasing functions on $\cP$ 
 which depend only on the distance from the origin. 
 Then we can identify $\cS(\cP)^{\rm sym}$ with $\cS(\R)^{\rm even}$
 the set of rapidly decreasing even functions on $\R$. 

 On the other hand, Hilbert transform $\vp\rightarrow \cH\vp$ is involutive, 
 i.e. $\cH^2=\text{id}$. Moreover, $\cH$ exchanges odd functions and 
 even functions, or real valued functions and pure imaginary valued functions. 
 Getting together, the statement holds. 
\end{proof}

\begin{proof}[Proof of \ref{Thm:twistor_correspondence}]
The one-to-one correspondence is already given in Lemma \ref{Lem:h=>f}, 
\ref{Lem:f=>h} and \ref{Lem:one-to-one}. 
So we construct the double fibration 
$\widetilde{Gr}_2(\R^4)\leftarrow\hat{\cZ}\rightarrow\CP^3$ and 
we check the conditions (1), (2), (3) and (4). 
As explained in the proof of Proposition \ref{Prop:holo_disks}, 
the family of holomorphic disks in $(\CP^3,\hat{P})$ is parametrized by 
$\widetilde{Gr}_2(\R^4)$ for each $h\in i\cS(\R)^{\rm odd}$. 
Let $D_x$ be the holomorphic disk in $(\CP^3,P)$ 
which corresponds to $x\in\widetilde{Gr}_2(\R^4)$, and we put
$$\hat{\cZ}=\left\{(x,y)\in \widetilde{Gr}_2(\R^4)\times \CP^3 : 
    y\in D_x\right\}.$$
We define $\wp : \hat{\cZ}\rightarrow \widetilde{Gr}_2(\R^4)$ and 
$\Psi : \hat{\cZ}\rightarrow\CP^3$ as the projections. 
Respecting the natural embedding $\cZ\rightarrow\hat{\cZ}$, 
$\wp$ and $\Psi$ are natural extension of $\tilde{\Gp}$ and 
$\Phi_\C$ respectively. 
From Proposition \ref{Prop:holo_disks} and its proof, $\wp$ is a disk bundle 
with fiberwise complex structure, and $\Psi$ is a continuous surjection, 
so the condition (1) holds. 
The conditions (2) and (3) follow from Proposition \ref{Prop:holo_disks}, 
Lemma \ref{Lem:h=>f} and Lemma \ref{Lem:f=>h}. 

For the last, the condition (4) follows from the continuity 
of the family of holomorphic disks. 
Actually if $y\neq y_0$ then 
$\wp(\Psi^{-1}(y))$ is the closure of the $\b$-surface 
$\tilde{\Gp}(\Phi_\C^{-1}(y))\subset TS^2$, 
and if $y=y_0$ then $\wp(\Psi^{-1}(y_0))=S^2_\infty$ is the 
singular $\b$-surface. 
\end{proof}

From the condition (3) of Theorem \ref{Thm:twistor_correspondence}, 
$\wp$ induces a continuous map 
$\varpi : (\CP^3-\hat{P}) \rightarrow\widetilde{Gr}_2(\R^4)$ 
which is smooth on $\cO(1)-\cO(1)|_{\RP^2}$. 
Next proposition says that $\varpi$ is not differentiable 
on $\cO(1)|_{\RP^2}-P$ when the singularity exists. 

\begin{Prop}
 $\varpi$ is differentiable if and only if $h(\mu)\equiv 0$. 
\end{Prop}
\begin{proof}
 Let $V=\{[1:w_1:w_2:w_3]\in\CP^3\}$ be an affine open set. 
 Then $V\cap \hat{P}$ is given by 
 $\{(w_1,w_2,w_3) : \Imag w_1=\Imag w_2=\Imag w_3-\Imag H=0\}$, 
 where $H$ is a function on $\{(w_1,w_2) : \Imag w_1\neq 0\}$ given by 
 $$ H=H^\pm\( \frac{\Imag \bar{w}_1w_2}{\Imag w_1}, 
    -\frac{\Imag w_2}{\Imag w_1}, w_1 \) $$ 
 by using $H^\pm$ defined in (\ref{Phi_TD+}) and so on, and 
 $\Imag H$ extends continuously on $\{\Imag w_1=0\}$ from the definition. 
 Let $U$ be an open neighborhood of $\widetilde{Gr}_2(\R^4)$ 
 given by (\ref{Prop:singularity_of_metric}). 
 Then the image $\varpi(V\setminus\hat{P})$ is contained in $U$ 
 and $\varpi|_{V\setminus\hat{P}}$ is described by 
 $$\begin{aligned} 
     u_1=\frac{\Imag \bar{w}_1(w_3-H)}{\Imag(w_3-H)}, \quad
     u_2=-\frac{\Imag w_1}{\Imag(w_3-H)}, \\
	 u_3=-\frac{\Imag w_2}{\Imag(w_3-H)}, \quad
	 u_4=\frac{\Imag \bar{w}_2(w_3-H)}{\Imag(w_3-H)}. 
	 \end{aligned} $$
Remark that these equations are defined on 
$\{\Imag w_1\neq 0\}\cap (V\setminus \hat{P})$, 
however these are continuously extended on $V\setminus \hat{P}$. 

Now we prove that if $\varpi$ is differentiable, 
then, for all $A>0$, $h(\mu)=0$ for any $\mu\in[-A,A]$. 
Let $t\in(-\e,\e)$ be a parameter on a small interval, and we fix $s\in\R$. 
We take a smooth curve in $V\setminus\hat{P}$ given by 
$$(w_1,w_2,w_3)=(s+it,A(-s+it),c(t)), $$ 
where $c(t)$ is a smooth function. Then 
$$H=H(t)=H^\pm(0,-A,s+it)$$ 
is defined on $t\neq 0$ and $\Imag H(t)$ is defined on $t\in(-\e,\e)$. 
We have $\Imag c(t)-\Imag H(t)\neq 0$ for all $t\in(-\e,\e)$, and 
$$u_4=u_4(t) = -\frac{\Real (c(t)-H(t))}{\Imag (c(t)-H(t))}\, t -s. $$ 
Hence
$$ \lim_{t\rightarrow 0} \frac{du_4}{dt} = 
 -\lim_{t\rightarrow 0} \frac{\Real (c(t)-H(t))}{\Imag (c(t)-H(t))}. $$
Since $\Imag H(t)$ is continuous, $\Real H(t)$ is also 
continuous if $\varpi$ is differentiable. 
By definition, we have 
$$\begin{aligned}
 H(t) &= H^+(0,-A,s+it)=F\(\frac{A}{2},\w\)
     =\frac{4i}{1-\w} \sum_{l=0}^\infty h_{\frac{A}{2},l}\,\w^{l+1} 
	  & \text{on} \quad t>0, \\
 H(t) &= H^-(0,-A,s+it)=-\overline{F\(\frac{A}{2},\bar{w}^{-1}\)}
     =-\frac{4i}{1-\w} \sum_{l=0}^\infty \overline{h_{\frac{A}{2},l}}\,\w^{-l} 
	  & \text{on} \quad t<0, \end{aligned}$$
where 
$\w=\w(t)=\frac{\textstyle s+it-i}{\textstyle s+it+i}. $
If we evaluate $a=\frac{A}{2}$ to the expansion of $h$ 
given in (\ref{expand_h}), we have 
\begin{equation} \label{eq:h(Acos)}
 h\(A\cos\frac{\c}{2}\)=\sum_{l=0}^\infty 
 \left\{ h_{\frac{A}{2},l}e^{i\c(l+\frac{1}{2})}- 
 \overline{h_{\frac{A}{2},l}}e^{-i\c(l+\frac{1}{2})} \right\}. 
\end{equation}
Then we have 
$$ \begin{aligned}
  \, & \lim_{t\searrow 0} \Imag H(t)= \lim_{t\nearrow 0} \Imag H(t) =
    \frac{i}{\sin\frac{\c}{2}}\, h\(A\cos\frac{\c}{2}\), \\
  \, & \lim_{t\searrow 0} \Real H(t)= -\lim_{t\nearrow 0} \Real H(t) =
    -\frac{1}{\sin\frac{\c}{2}} \sum_{l=0}^\infty
	\left\{ h_{\frac{A}{2},l}e^{i\c(l+\frac{1}{2})}+ 
  \overline{h_{\frac{A}{2},l}}e^{-i\c(l+\frac{1}{2})} \right\},  
  \end{aligned} $$
where we define $\c=\c(s)\in (0,2\pi)$ so that $e^{i\c}=\frac{s-i}{s+i}=\w(0)$. 
The first equation says that $\Imag H(t)$ is indeed continuous. 
The summand in the right hand side of the second equation defines
a Fourier expansion for some $L^2$-function for $e^{\frac{i\c}{2}}$ which 
is zero almost everywhere if $\Real H(t)$ is continuous for every $s\in\R$. 
Hence $h_{\frac{A}{2},l}=0$ for all $l$, and we have 
$h(A\cos\frac{\c}{2})=0$ for all $\c$ from (\ref{eq:h(Acos)}). 
Therefore, for all $A>0$, 
$h(\mu)=0$ for any $\mu\in[-A,A]$ if $\varpi$ is differentiable. 
\end{proof}

\section{Appendix : Radon transform} \label{Appendix}
Here is a review of the Radon transform over $\R^2$ (cf.\cite{bib:Helgason}). 
Let $\tilde{\cP}$ be the set of oriented lines in $\R^2$. 
Then $\tilde{\cP}$ is diffeomorphic to $S^1\times \R$, 
where the correspondence 
$(\R/2\pi\Z)\times\R\ni(\s,\mu) \mapsto \x\in\tilde{\cP}$ 
is given by $\x=\{t(\cos\s,\sin\s)+\mu(-\sin\s,\cos\s)\}$ 
using parameter $t\in\R$. 
Let $\cP$ be the set of unoriented lines in $\R^2$, then 
$\cP\cong(S^1\times\R)/\Z_2$ where the equivalence is 
$(\s,\mu)\sim(\s+\pi,-\mu)$. 
Let $f$ be a rapidly decreasing complex valued function over $\R^2$, 
then the {\it Radon transform} $\hat{f}$ of $f$ is defined by 
\vspace{-1mm}
\begin{equation}
  \hat{f}(\s,\mu)=\int_{-\infty}^\infty 
  f(t\cos\s-\mu\sin\s,t\sin\s+\mu\cos\s)dt. 
\end{equation}
Then $\hat{f}$ is rapidly decreasing function on $\cP$. 
The {\it dual transform} $\vp^\vee$ of rapidly decreasing function 
$\vp$ on $\cP$ is defined by
\vspace{-1mm}
\begin{equation}
 \vp^\vee (x)= \frac{1}{2\pi} \int_{S^1} \vp(\s,\mu(x,\s))d\s, 
\end{equation}
where $\mu(x,\mu)=-x_1\sin\s+x_2\cos\s$. 
Remark that $(\s,\mu(x,\s))$ runs all the oriented lines through $x\in\R^2$. 
The inversion formula is given by
\vspace{-1mm}
\begin{equation} \label{inversion_formula}
 f= \frac{1}{2i} \(\pd{\mu} {\mathcal H}_\mu\hat{f} \)^\vee,  
\end{equation}
\vspace{-1mm}
where ${\mathcal H}_\mu$ is the {\it Hilbert transform} defined by 
\vspace{-1mm}
\begin{equation} \label{Hilbert_transf}
  ({\mathcal H}_\mu \vp)(\s,\mu)=\frac{i}{\pi}\ \pv \int_{-\infty}^\infty 
  \frac{\vp(\s,\nu)}{\nu-\mu}d\nu. 
\end{equation}

\vspace{1ex}
\noindent
{\large \bf Acknowledgments} \\
The author most gratefully thanks his supervisor Mikio Furuta, 
as well as Hiroyuki Kamada and Yasuo Matsushita for 
comments and conversations. 

\renewcommand{\bibname}{ŽQl•¶Œ£}

\end{document}